\documentclass[final,1p,12pt]{elsarticle}

\usepackage{hyperref}
\usepackage{amssymb}
\usepackage{amsmath}
\usepackage{dirtytalk}
\usepackage{bm}
\usepackage{todonotes}
\usepackage{dutchcal}
\usepackage{nicefrac}
\usepackage{mathtools}
\usepackage{listofitems}
\usepackage[algo2e,ruled,vlined]{algorithm2e}
\tikzstyle{mynode}=[thick,draw=blue,fill=blue!20,circle,minimum size=30]
\usepackage{cleveref}

\newcommand{\cor}[1]{#1}
\newenvironment{core}
    {}
    {}
\newcommand{\corr}[1]{#1}
\newcommand{\corrr}[1]{#1}

\usepackage{lineno}

\journal{CMAME}
\begin{document}
\begin{frontmatter}



\title{Multilevel domain decomposition-based architectures for physics-informed neural networks}


\author[a1]{Victorita Dolean} 
\affiliation[a1]{organization={Centre for Analysis Scientific Computing and Applications (CASA), Eindhoven University of Technology},
            addressline={PO Box 513, 5600 MB Eindhoven}, 
            country={The Netherlands}}

\author[a2]{Alexander Heinlein} 
\affiliation[a2]{organization={Delft Institute of Applied Mathematics, Delft University of Technology},
            addressline={Mekelweg 4, 2628 CD Delft}, 
            country={The Netherlands}}

\author[a3]{Siddhartha Mishra} 
\author[a3]{Ben Moseley} 
\affiliation[a3]{organization={Seminar for Applied Mathematics (SAM) D-MATH, and ETH AI Center},
            addressline={ETH Zurich}, 
            country={Switzerland}}

\begin{abstract}
Physics-informed neural networks (PINNs) are a powerful approach for solving problems involving differential equations, yet they often struggle to solve problems with high frequency and/or multi-scale solutions. Finite basis physics-informed neural networks (FBPINNs) improve the performance of PINNs in this regime by combining them with an overlapping domain decomposition approach. In this work, FBPINNs are extended by adding multiple levels of domain decompositions to their solution ansatz, inspired by classical multilevel Schwarz domain decomposition methods (DDMs). Analogous to typical tests for classical DDMs, we assess how the accuracy of PINNs, FBPINNs and multilevel FBPINNs scale with respect to computational effort and solution complexity by carrying out strong and weak scaling tests. Our numerical results show that the proposed multilevel FBPINNs consistently and significantly outperform PINNs across a range of problems with high frequency and multi-scale solutions. Furthermore, as expected in classical DDMs, we show that multilevel FBPINNs improve the accuracy of FBPINNs when using large numbers of subdomains by aiding global communication between subdomains.
\end{abstract}


\begin{highlights}
\item We propose a method for solving differential equations.
\item Our method combines PINNs with multilevel domain decomposition.
\item Our approach significantly outperforms PINNs when solving multi-scale problems.
\item Multilevel modeling improves accuracy by aiding communication between subdomains.
\end{highlights}


\end{frontmatter}


\section{Introduction} \label{sec:introduction}




Scientific machine learning (SciML) \cite{Baker2019, Willard2022, Cuomo2022a, Arridge2019, Moseley2022} is an emerging and rapidly growing field of research. The central goal of SciML is to provide accurate, efficient, and robust tools for carrying out scientific research by tightly combining scientific understanding with machine learning (ML). The field has provided many such tools which have enhanced traditional approaches, from accelerating simulation algorithms to discovering new scientific phenomena.


One popular SciML approach \corrr{is} physics-informed neural networks (PINNs) \cite{Lagaris1998, Raissi2019}. PINNs solve forward and inverse problems related to differential equations by using a neural network to directly approximate the solution to the differential equation. They are trained by using a loss function which minimizes the residual of the differential equation over a set of collocation points. The initial concepts behind PINNs were introduced by~\cite{Lagaris1998} and others, and later re-implemented and extended in~\cite{Raissi2019}. One of the advantages of PINNs over traditional methods for solving differential equations such as finite difference (FD) and finite element methods (FEM) is that they provide a mesh-free approach, paving the way for the application of problems with complex geometry or in very high spatial dimensions; cf.~\cite{MM3}. Furthermore, they can easily be extended to solve inverse problems by incorporating observational data.

Since their invention, PINNs have been employed across a wide range of domains \cite{Cuomo2022a, Karniadakis2021}. For example, they have been used to solve forward and inverse problems in geophysics~\cite{Moseley2020c}, fluid dynamics~\cite{Jin2021,Cai2021,Raissi2020}, and optics~\cite{Chen2020}. Many extensions of PINNs have also been proposed. For example, PINNs have been extended to carry out uncertainty quantification~\cite{Yang2021}, learn fast surrogate models~\cite{Wang2021e,Zhu2019b}, and carry out equation discovery~\cite{Chen2021b}.


However, PINNs suffer from a number of limitations. One is that, compared to traditional methods, their convergence properties are poorly understood, although some work \corrr{has} started to explore this \cite{Mishra2022, Shin2020, Wang2022}. Another limitation is that, compared to traditional methods, the computational cost of training PINNs is relatively high, especially when they are only used for forward modeling \cite{Karniadakis2021}. Finally, a major limitation of PINNs is that they often struggle to solve problems with high frequency and/or multi-scale solutions~\cite{Moseley2023,Wang2021d}. Typically, as higher frequencies and multi-scale features are added to the solution, the accuracy of PINNs usually rapidly reduces and their computational cost rapidly increases in a super-linear fashion \cite{Moseley2023}.


There are multiple reasons for this behavior. One is the spectral bias of neural networks, which is the well-studied property that neural networks tend to learn high frequencies much slower than low frequencies~\cite{Xu2020e,Rahaman2018,Basri2019,Cao2019}. Another is that, as higher frequencies and more multi-scale features are added, more collocation points and a larger neural network with significantly more free parameters are typically required to accurately approximate the solution. This creates a significantly more complex optimization problem when training the PINN.


Recently, \cite{Moseley2023} proposed finite basis physics-informed neural networks (FBPINNs), which aim to improve the performance of PINNs in this regime by using an overlapping domain decomposition (DD) approach. In particular, instead of using a single neural network to approximate the solution to the differential equation, many smaller neural networks were placed in overlapping subdomains and summed together to represent the solution. On the one hand, FBPINNs can be seen as a DD-based network architecture for PINNs. On the other hand, by taking this \say{divide and conquer} approach, the global PINN optimization problem is transformed into many smaller local optimization problems, which are coupled implicitly due to the overlap of the subdomains and their globally defined loss function. The results in~\cite{Moseley2023} show that this approach significantly improves the accuracy and reduces the training cost of PINNs when solving differential equations with high frequency and multi-scale solutions.

In this work, we significantly extend FBPINNs by incorporating \emph{multilevel modeling} into their design. In particular, instead of using a single DD in their solution ansatz, we add multiple levels of overlapping DDs. This idea is inspired by classical DDMs, where coarse levels are required for numerical scalability when using large numbers of subdomains. Furthermore, to assess the performance of multilevel FBPINNs, we define strong and weak scaling tests for measuring how the accuracy of PINNs and FBPINNs scale with computational effort and solution complexity, analogous to the strong and weak scaling tests commonly used in classical DDMs.

Given these extensions, the performance of PINNs, (one-level) FBPINNs, and multilevel FBPINNs across a range of high frequency and multi-scale problems is investigated. \cor{We also compare multilevel FBPINNs to PINNs with Fourier input features \cite{Tancik2020, Wang2021d} and self-adaptive PINNs (SA-PINNs) \cite{McClenny2023}, which have both been shown to improve the accuracy of PINNs when solving multi-scale problems}. Across these tests, we find that multilevel FBPINNs significantly outperform PINNs in terms of accuracy and computational cost. Furthermore, as expected in classical DDMs, we show that multilevel FBPINNs improve the accuracy of FBPINNs when using large numbers of subdomains by aiding global communication between subdomains.

The remainder of this work is structured as follows. In~\cref{sec:related_work} we discuss related work on combining ML, PINNs, and DD, and in~\cref{sec:nns,sec:pinns} we give a brief overview of neural networks and PINNs. Then we define FBPINNs and extend them to multilevel FBPINNs in~\cref{sec:methods}. Our strong and weak scaling tests and corresponding numerical results on the performance of PINNs, FBPINNs, and multilevel FBPINNs across a range of high frequency and multi-scale problems are discussed in~\cref{sec:results}. Finally, in~\cref{sec:discussion}, we discuss the implications and limitations of our work and further research directions.

\subsection{Related work} \label{sec:related_work}

In general, the idea of combing ML with classical DDMs is not new; for early works on using ML to predict the geometrical location of constraints in adaptive finite element tearing and interconnecting (FETI) and balancing DD by constraints (BDDC) methods; see~\cite{Heinlein:2019:MLA}. An overview of the first attempts on combining DD and ML can be found in~\cite{Heinlein:2021:CML2}\cor{, a more recent overview is given in~\cite{klawonn2023machine}}. 

For specifically combining PINNs with DD, some of the first methods in this area were the deep domain decomposition method (D3M)~\cite{li2019d3m}, the deep-learning-based domain decomposition method (DeepDDM) \cite{li2020deep, li_deep_2023b}, and its two-level variant \cite{mercier:2021:ACS}, which use PINNs to solve local problems and overlapping Schwarz steps to iteratively connect them based on Lions' parallel Schwarz algorithm~\cite{lions_schwarz_1988}. At the same time, a series of other extensions, like \cor{conservative physics-informed neural network (cPINN) and extended physics-informed neural networks (XPINNs)}~\cite{Jagtap2020b} were proposed, which similarly divide the domain and use PINNs to solve each local problem; here, typically a non-overlapping DD is used. \cor{A detailed comparison of these methods to FBPINNs is given in \cref{sec:fbpinn-dd-compare}.}

In~\cite{lee_partition_2021}, partition of unity functions\corrr{,} similar to the window functions used in the FBPINN method, are learned. However, this is done in a pure function approximation setting rather than in the solution of PDE-based problems with PINNs.

\subsection{Neural networks} \label{sec:nns}

We first provide a basic definition of a neural network. For the purpose of this work, we simply consider a neural network to be a mathematical function with some learnable parameters. More precisely, the network is defined as $u(\mathbf{x}, \bm{\theta}): \mathbb{R}^{d_x} \times \mathbb{R}^{d_\theta} \rightarrow \mathbb{R}^{d_u}$, where $\mathbf{x}$ are some inputs to the network, $\bm{\theta}$ are a set of learnable parameters, and $d_x$, $d_\theta$, and $d_u$ are the dimensionality of the network's inputs, parameters, and outputs. In a traditional supervised learning setting, learning typically consists of fitting the network function to some training data containing example inputs and outputs, by minimizing a loss function with respect to $\bm{\theta}$ which penalizes the difference between the network's outputs and the training data.

The exact form of the network function is determined by the neural network's architecture. In this work, we solely use feedforward fully connected networks (FCNs) \cite{Goodfellow2016}. In this case, the network function is given by
\begin{equation} \label{eq:fcn}
u(\mathbf{x}, \bm{\theta}) = f_n \circ ... \circ f_i \circ ... \circ f_1(\mathbf{x}, \bm{\theta})
\end{equation}
where now $\mathbf{x} \in \mathbb{R}^{d_0}$ is the input to the FCN,
$u \in \mathbb{R}^{d_n}$ is the output of the FCN, 
$n$ is the number of layers (depth) of the FCN, 
and $f_i(\mathbf{x}, \bm{\theta})=\sigma_i(W_i\mathbf{x}+\mathbf{b}_i)$
where
$\bm{\theta}_i = (W_i, \mathbf{b}_i)$, 
$W_i \in \mathbb{R}^{d_i \times d_{i-1}}$ are known as weight matrices, $\mathbf{b}_i \in \mathbb{R}^{d_i}$ are known as bias vectors,
$\sigma_i$ are element-wise activation functions commonly chosen as rectified linear unit (ReLU), hyperbolic tangent, or identity functions,
and $\bm{\theta} = (\bm{\theta}_1, ..., \bm{\theta}_i , ..., \bm{\theta}_n)$ are the set of learnable parameters of the network. Note that only the nonlinear activation functions $\sigma_i$ facilitate nonlinearity of the network function.

\subsection{Physics-informed neural networks} \label{sec:pinns}

Physics-informed neural networks (PINNs) \cite{Lagaris1998, Raissi2019} use neural networks to solve problems related to differential equations. In particular, PINNs focus on solving boundary value problems of the form
\begin{equation} \label{eq:bvp}
\begin{array}{rcll}
\mathcal{N}[u]({\mathbf{x}})    & = & f(\mathbf{x}), \text{ }   & \mathbf{x} \in \Omega \subset \mathbb{R}^{d},     \\
\mathcal{B}_k[u] ({\mathbf{x}}) & = & g_k(\mathbf{x}), \text{ } & \mathbf{x}  \in \Gamma_{k} \subset \partial\Omega
\end{array}
\end{equation}
where $\mathcal{N}[u](\mathbf{x}) $ is some differential operator, $u(\mathbf{x})$ is the solution, and $\mathcal{B}_k(\cdot)$ are a set of boundary conditions (BCs) which ensure uniqueness of the solution. For the sake of simplicity, we consider BCs in a broad sense; we do not explicitly distinguish between initial and boundary conditions, and the $\mathbf{x}$ variable can include time. \cref{eq:bvp} can describe many different differential equation problems, including linear and non-linear problems, time-dependent and time-independent problems, and those with irregular, higher-order, and cyclic boundary conditions.

To solve \cref{eq:bvp}, PINNs use a neural network to directly approximate the solution, i.e., $u(\mathbf{x}, \bm{\theta})\approx u(\mathbf{x})$. Note, for simplicity throughout this work, we use the same notation for the true solution and the neural network. It is important to note that PINNs provide a functional approximation to the solution, and not a discretized solution similar to that provided by traditional methods such as finite difference methods, and as such PINNs are a mesh-free approach for solving differential equations. Following the approach proposed by \cite{Raissi2019}, the following loss function is minimized to train the PINN,
\begin{equation} \label{eq:pinn_loss_soft}
\mathcal{L}(\bm{\theta}) 
= 
\frac{\lambda_{I}}{N_I}
\sum_{i=1}^{N_I}
\big(
\underbrace{
    \mathcal{N}[u](\mathbf{x}_{i},\bm{\theta}) - f(\mathbf{x}_{i})
}_{\text{PDE residual}} 
\big)^{2}
+ 
\sum_{k=1}^{N_k} 
\frac{\lambda_{B}^{k}}{N_{B}^{k}} \sum_{i=1}^{N_{B}^{k}}
\big(
\underbrace{
    \mathcal{B}_k[u](\mathbf{x}_{i}^{k},\bm{\theta})  - g_k(\mathbf{x}_{i}^{k})
}_{\text{BC residual}}
\big)^2. 
\end{equation}
where $\{\mathbf{x}_i\}_{i=1}^{N_I}$ is a set of collocation points sampled in the interior of the domain, $\{\mathbf{x}_{j}^{k}\}_{j=1}^{N_{B}^{k}}$ is a set of points sampled along each boundary condition, and $\lambda_I$ and $\lambda_{B}^{k}$ are well-chosen scalar weights that ensure the terms in the loss function are well balanced. Intuitively, one can see that by minimizing the PDE residual, the method tries to ensure that the solution learned by the network obeys the underlying PDE, and by minimizing the BC residual, the method tries to ensure that the learned solution is unique by matching it to the BCs. Importantly, a sufficient number of collocation and boundary points must be chosen such that the PINN is able to learn a consistent solution across the domain.

Iterative schemes are typically used to optimize this loss function. Usually, variants of the gradient descent (GD) method, such as the Adam optimizer~\cite{kingma_adam_2017}, or quasi-Newton methods, such as the limited-memory Broyden--Fletcher--Goldfarb--Shanno (L-BFGS) algorithm~\cite{liu_limited_1989} are employed. These methods require the computation of the gradient of the loss function with respect to the network parameters, which can \cor{be} computed easily and efficiently using automatic differentiation~\cite{kelley_gradient_1960} provided in modern deep learning libraries \cite{abadi_tensorflow_2016, paszke_pytorch_2019, jax}. Note that gradients of the network output with respect to its inputs are also typically required to evaluate the PDE residual in the loss function, and can similarly be obtained and further differentiated through to update the network's parameters using automatic differentiation.

\subsubsection{Hard constrained PINNs} \label{sec:hard_pinns}

A downside of training PINNs with the loss function given by \cref{eq:pinn_loss_soft} is that the BCs are {\it softly} enforced. This means the learned solution may deviate from the BCs because the BC term may not be fully minimized. Furthermore, it can be challenging to balance the different objectives of the PDE and BC terms in the loss function, which can lead to poor convergence and solution accuracy \cite{Wang2022, Sun2019a}. An alternative approach, as originally proposed by \cite{Lagaris1998}, is to enforce BCs in a {\it hard} fashion by using the neural network as part of a solution ansatz. More precisely, the solution to the differential equation is instead approximated by $[\mathcal{C} u](\mathbf{x},\bm{\theta}) \approx u(\mathbf{x})$ where $\mathcal{C}$ is an appropriately selected constraining operator which analytically enforces the BCs~\cite{Moseley2023, Leake2020}. 

To give a simple example, suppose we want to enforce $u(x=0)=0$ when solving a one-dimensional ordinary differential equation (ODE). The constraining operator and solution ansatz could be chosen as $[\mathcal{C} u] (x,\bm{\theta}) = \tanh(x) u(x,\bm{\theta})\approx u(\mathbf{x})$.
The rationale behind this is that the function $\tanh(x) $ is zero at $0$, forcing the BC to always be obeyed, but non-zero away from $0$, allowing the network to learn the solution away from the BC.

In this approach, the BCs are always satisfied and therefore the BC term in the loss function \cref{eq:pinn_loss_soft} can be removed, meaning that the PINN can be trained using the simpler unconstrained loss function,
\begin{equation} \label{eq:pinn_loss_hard}
\mathcal{L}(\bm{\theta}) =\frac{1}{N} \sum_{i=1}^{N} (\mathcal{N}[\mathcal{C} u](\mathbf{x}_{i},\bm{\theta}) - f(\mathbf{x}_{i}))^{2}.
\end{equation}
where $\{\mathbf{x}_i\}_{i=1}^{N}$ is a set of collocation points sampled in the interior of the domain. Note that, in general, there is no unique way of choosing the constraining operator, and the definition of a suitable constraining operator for complex geometries and/or complex BCs may be difficult or sometimes even impossible, i.e., this strategy is problem dependent; in this case, one may resort to the soft enforcement of boundary conditions~\cref{eq:pinn_loss_soft} instead.

\section{Methods} \label{sec:methods}

In this section, we define FBPINNs (\cref{sec:fbpinns}) and extend them to multilevel FBPINNs (\cref{sec:multilevel-fbpinns}). We also discuss the similarities and differences of FBPINNs and multilevel FBPINNs to classical DDMs (\cref{sec:classical-comparison}).

\subsection{Finite basis physics-informed neural networks} \label{sec:fbpinns}

\begin{figure}[!t] 
\centering
\includegraphics[width=0.6\textwidth]{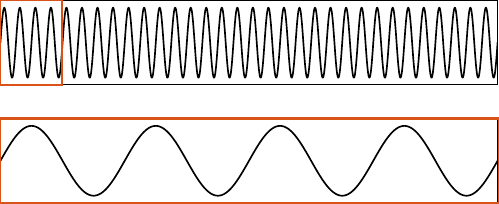}
\caption{Scaling high frequency problems to low frequency problems using domain decomposition. FBPINNs decompose the domain into many subdomains, and use neural networks within each subdomain to learn the local solution. The input coordinates to each network are normalized to the range [-1,1] over their individual subdomains. When solving problems with high frequency solutions, this effectively scales each local problem from a high frequency problem to a lower frequency problem, and helps reduce the network's spectral bias.
\label{fig:high_frequency}}
\end{figure}

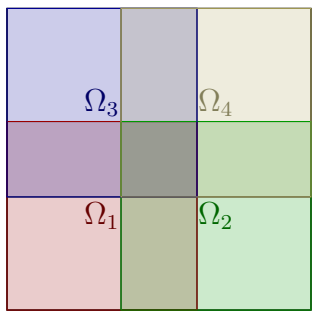
\begin{figure}[!t]
\centering
\begin{tikzpicture}
    \draw (0,0) -- +(4,0) -- +(4,4) -- +(0,4) -- cycle;
    
    \draw[fill=red!60!black,fill opacity=0.2] (0,0) -- +(2.5,0) -- +(2.5,2.5) -- +(0,2.5) -- cycle;
    \draw[fill=green!60!black,fill opacity=0.2] (1.5,0) -- +(2.5,0) -- +(2.5,2.5) -- +(0,2.5) -- cycle;
    \draw[fill=blue!60!black,fill opacity=0.2] (0,1.5) -- +(2.5,0) -- +(2.5,2.5) -- +(0,2.5) -- cycle;
    \draw[fill=yellow!60!black,fill opacity=0.2] (1.5,1.5) -- +(2.5,0) -- +(2.5,2.5) -- +(0,2.5) -- cycle;

    \draw[red!60!black] (0,0) -- +(2.5,0) -- +(2.5,2.5) -- +(0,2.5) -- cycle; 
    \node[red!40!black] (omega1) at (1.25,1.25) {$\Omega_1$};
    \draw[green!60!black] (1.5,0) -- +(2.5,0) -- +(2.5,2.5) -- +(0,2.5) -- cycle;
    \node[green!40!black] (omega1) at (2.75,1.25) {$\Omega_2$};
    \draw[blue!60!black] (0,1.5) -- +(2.5,0) -- +(2.5,2.5) -- +(0,2.5) -- cycle;
    \node[blue!40!black] (omega1) at (1.25,2.75) {$\Omega_3$};
    \draw[yellow!60!black] (1.5,1.5) -- +(2.5,0) -- +(2.5,2.5) -- +(0,2.5) -- cycle;
    \node[yellow!40!black] (omega1) at (2.75,2.75) {$\Omega_4$};
\end{tikzpicture}
\caption{Plot of a square domain $\Omega$ decomposed into four overlapping subdomains, using a uniform rectangular decomposition.
\label{fig:overlapping_dd}
}
\end{figure}

As discussed in \cref{sec:introduction}, a major challenge when training PINNs is that, when higher frequencies and multi-scale features are added to the solution, the accuracy of PINNs usually rapidly reduces and their computational cost rapidly increases in a super-linear fashion \cite{Moseley2023, Wang2021d}.

In the FBPINN approach \cite{Moseley2023}, instead of using a single neural network to represent the solution, many smaller neural networks are confined in overlapping subdomains and summed together to represent the solution. By taking this \say{divide and conquer} approach, the global PINN optimization problem is transformed into many smaller coupled local optimization problems. 

Furthermore, FBPINNs ensure that the inputs to each subdomain network are normalized over their individual subdomain. When solving problems with high frequency solutions, this effectively scales each local problem from a high frequency problem to a lower frequency problem, and helps limit the effect of spectral bias; \cref{fig:high_frequency} explains this effect further.

\subsubsection{Mathematical definition}

We now provide a mathematical definition of FBPINNs. First, the global solution domain $\Omega$ is decomposed into $J$ overlapping subdomains $\left\lbrace \Omega_j \right\rbrace_{j=1}^J$; cf. \cref{fig:overlapping_dd}. Then, for each subdomain $\Omega_j$, a space of network functions is defined,
$$
\mathcal{V}_j =
\left\lbrace \mathcal{v}_j(\mathbf{x}, \bm{\theta}_j) \mid \mathbf{x} \in \Omega_j, \bm{\theta}_j \in \Theta_j \right\rbrace,
$$
where $\mathcal{v}_j(\mathbf{x}, \bm{\theta}_j)$ is a neural network placed in each subdomain and $\Theta_j = \mathbb{R}^{K_j}$ is the linear space of all possible network parameters. Here, $K_j$ is the number of local network parameters which is determined by the network architecture.

Next, each subdomain network is confined to its subdomain by multiplying each network with a window function $\omega_j(\mathbf{x})$, where ${\rm supp} \big( \omega_j \big) \subset \Omega_j$. Note the neural network functions used in $\mathcal{V}_j$ generally can have global support, and the window functions are used to restrict them to their individual subdomains. Furthermore, we impose that the window functions form a partition of unity, i.e.,
$$
\sum_{j=1}^J \omega_j \equiv 1
\quad
\text{on } \Omega.
$$

Given the space of network functions and the window functions, we define a global space decomposition given by $\mathcal{V}$ as
$$
\mathcal{V}
=
\sum_{j=1}^{J} \omega_j \mathcal{V}_j.
$$
This space decomposition allows for decomposing any given function $u \in \mathcal{V}$ as follows
\begin{equation} \label{eq:fbpinn_network_architecture}
u
=
\sum_{j=1}^{J} \omega_j \mathcal{v}_j
\quad
\text{or}
\quad
u(\mathbf{x}, \bm{\theta})
=
\sum_{j=1}^{J} \omega_j \mathcal{v}_j(\mathbf{x}, \bm{\theta}_j),
\end{equation}
respectively.

FBPINNs solve the boundary value problem~\cref{eq:bvp} by using equation~\cref{eq:fbpinn_network_architecture} to approximate the solution, and we refer to \cref{eq:fbpinn_network_architecture} as the FBPINN solution. From a PINN perspective, the FBPINN solution can simply be thought of as a specific type of neural network architecture for the PINN which sums together many locally-confined networks to generate the output solution.

The same scheme for training PINNs is used to train the FBPINN. More specifically, the FBPINN solution \cref{eq:fbpinn_network_architecture} is substituted into the PINN loss function \cref{eq:pinn_loss_soft} and the same iterative optimization scheme is used to learn the parameters $\left\lbrace \bm{\theta}_j \right\rbrace_{j=1}^J$ of each subdomain network. FBPINNs can also be trained with hard BCs by using the same constraining operator approach described in \cref{sec:hard_pinns}. In particular, substituting the FBPINN solution \cref{eq:fbpinn_network_architecture} into the hard-constrained loss function \cref{eq:pinn_loss_hard} yields the loss function
\begin{equation} \label{eq:fbpinn_loss}
\mathcal{L}(\bm{\theta}) =\frac{1}{N} \sum_{i=1}^{N} (\mathcal{N}[\mathcal{C} \sum\limits_{j=1}^J \omega_j \mathcal{v}_j(\mathbf{x}_{i},\bm{\theta}_j)] - f(\mathbf{x}_{i}))^{2}.
\end{equation}

\subsubsection{Computational efficiency of FBPINNs versus PINNs} \label{sec:fbpinn-efficiency}
Assuming 
the same size network in each subdomain, naively computing the FBPINN solution \cref{eq:fbpinn_network_architecture} has a time complexity of $\mathcal{O}(N J \tilde S)$, where 
$\tilde S$ is the cost of computing the output of a single subdomain network for a single collocation point. This becomes very expensive as more subdomains are added ($J$ increases). However, because the output of each subdomain network is zero outside of the overlapping subdomain after applying the window function, only collocation points within the subdomain need to be included in the summations in \cref{eq:fbpinn_loss}. This reduces the computational cost to $\mathcal{O}(NC \tilde S)$, where $C$ is the average number of subdomains a collocation point belongs to. PINNs have a computational cost of $\mathcal{O}(NS)$, where $S$ is the cost of computing the output of the single global PINN network. Importantly, as the problem complexity increases, the size of the PINN network must typically be increased (increasing $S$), whilst for FBPINNs, we can typically keep the subdomain network size fixed, and increase $J$ instead -- thus, often $\tilde S \ll S$ and FBPINNs are often orders of magnitude more efficient than PINNs. More details on our efficient software implementation are provided in~\cref{sec:software}.

\begin{core}
\subsubsection{FBPINNs versus other methods for combining PINNs with DD} \label{sec:fbpinn-dd-compare}

Multiple other approaches exist which combine PINNs with DD; a recent overview can be found in~\cite{klawonn2023machine}. Similar to FBPINNs, XPINNs \cite{Jagtap2020b}  divide the domain into subdomains and use separate neural networks to solve each subdomain problem. However, XPINNs use a non-overlapping DD. A downside of this approach is that the global solution contains discontinuities at the subdomain interfaces and additional loss terms are required to enforce coupling between subdomain networks. In contrast, FBPINNs do not require additional loss terms and their solution is continuous across subdomain interfaces. XPINNs and FBPINNs are comparable in terms of their computational cost of training, and both are able to use irregular DDs with different types of subdomain neural networks. \cite{Dong2021} and \cite{Dwivedi2021} propose a similar approach to XPINNs, but they use extreme learning machines \cite{Huang2006} as subdomain networks, where only the parameters of the last layer of the network are learnt. This approach has the advantage of being much faster to train, but the capacity of the subdomain networks is limited, and the performance strongly depends on the initialization.

Other approaches attempt to learn suitable DDs for PINNs. For example, Gated-PINNs \cite{Stiller2020} use several neural networks, called experts, to propose a solution to a PDE, whilst a gating network is used to weight-average the expert solutions given a coordinate in the domain. Augmented-PINNs (APINNs)~\cite{Hu2022} build upon this strategy by introducing parameter sharing between experts to capture solution similarities between subdomains. These approaches are more flexible than FBPINNs in that they can adaptively learn DDs, but they are also significantly more computationally expensive to train as they require all experts to be evaluated at each input coordinate and therefore do not scale well with problem size; cf.~\cref{sec:fbpinn-efficiency}.

\end{core}

\subsection{Multilevel FBPINNs} \label{sec:multilevel-fbpinns}

\begin{figure}[!t]
\centering
\includegraphics[width=0.9\textwidth]{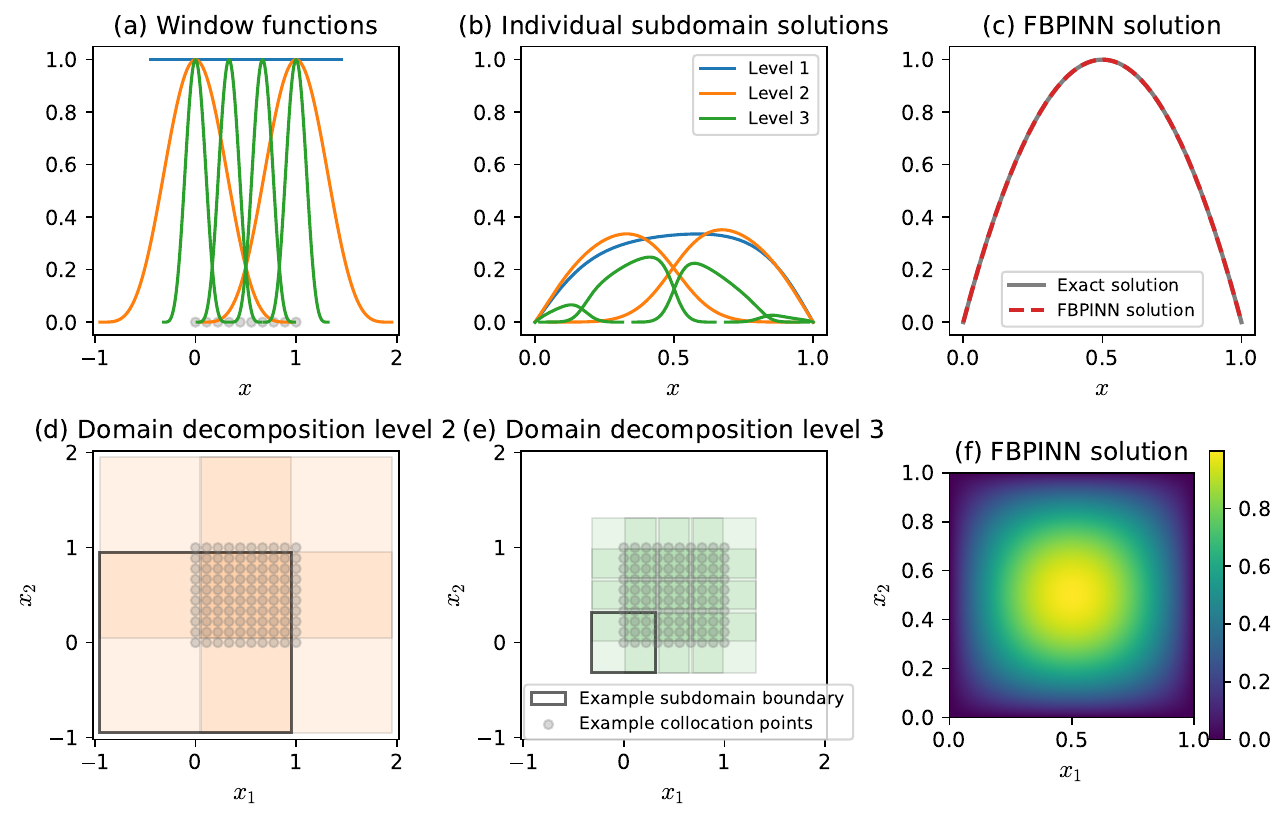}
\caption{Example of a multilevel FBPINN solving Laplace's equation in one and two dimensions. For the 1D problem, the multilevel FBPINN uses $L=3$ levels, where each level has $1$, $2$ and $4$ subdomains respectively. The window functions, $\hat \omega^{(l)}_j(x)$, used for each level are shown in (a), the individual solutions learned by each subdomain network are shown in (b), and the multilevel FBPINN solution is shown in (c). For the 2D problem, the multilevel FBPINN uses $L=3$ levels, where each level has $1\times1$, $2\times2$ and $4\times4$ subdomains respectively, using a uniform rectangular DD. The DDs for level $2$ and level $3$ are plotted in (d) and (e), and the multilevel FBPINN solution is shown in (f). Note the subdomain boundaries and window functions extend past the problem domain (in this case, $[0,1]^d$). Example collocation points used to train the multilevel FBPINN are plotted in (a), (d) and (e).
\label{fig:dd-collocation_points}}
\end{figure}

In this work we propose multilevel FBPINNs, which extend FBPINNs by adding multiple levels of DDs to their solution ansatz.  They are inspired by classical multilevel DD methods, where coarse levels are generally required for numerical scalability when using large numbers of subdomains, and multilevel approaches may significantly improve performance; see, for instance,~\cite{toselli_domain_2005,Dolean:AID:2015}. Our hypothesis is that adding multilevel modelling to FBPINNs similarly improves their performance. The generalization of FBPINNs to two levels was briefly discussed in~\cite{dolean:2022:FBP} and we fully introduce the concept here.

A multilevel FBPINN is defined as follows. First, we define $L$ levels of DDs, where each level, $l$, defines an overlapping DD of $\Omega$ with $J^{(l)}$ subdomains, i.e.,
$$
D^{(l)} = \left\lbrace \Omega_j^{(l)} \right\rbrace_{j=1}^{J^{(l)}},
$$
for $j = 1,\ldots,J^{(l)}$. Without loss of generality, let $J^{(1)} = 1$, that is, on the first level, we only have one subdomain $\Omega_j^{(1)} = \Omega$. Moreover, we let $J^{(1)} < J^{(2)} < \ldots < J^{(l)}$, meaning that the number of subdomains increases from one to the next level.

Next, we define spaces of network functions for each level,
$$
\mathcal{V}_j^{(l)} 
=
\left\lbrace \mathcal{v}_j^{(l)} (\mathbf{x}, \bm{\theta}_j^{(l)}) \mid \mathbf{x} \in \Omega_j^{(l)}, \bm{\theta}_j^{(l)}\in \Theta_j^{(l)} \right\rbrace,
\quad
j=1,\ldots,J^{(l)},\ l = 1,\ldots,L,
$$
as well as a partition of unity for each level using window functions, $\omega_j^{(l)}$, with
$$
{\rm supp} \big( \omega_j^{(l)} \big) \subset \Omega_j^{(l)}
\quad
\text{and}
\quad
\sum_{j=1}^{J^{(l)}} \omega_j^{(l)} \equiv 1
\quad
\text{ on } \Omega \quad \forall l.
$$
We can then define a global space decomposition,
$$
\mathcal{V}
=
\frac{1}{L} \sum_{l=1}^L \sum_{j=1}^{J^{(l)}} \omega_j^{(l)} \mathcal{V}_j^{(l)},
$$
and use this space decomposition to decompose any given function $u \in \mathcal{V}$ as follows,
\begin{equation} \label{eq:multilevel_fbpinn_network_architecture}
u
=
\frac{1}{L} \sum_{l=1}^L \sum_{j=1}^{J^{(l)}} \omega_j^{(l)} \mathcal{v}_j^{(l)}
\quad
\text{or}
\quad
u(\mathbf{x}, \bm{\theta})
=
\frac{1}{L} \sum_{l=1}^L \sum_{j=1}^{J^{(l)}} \omega_j^{(l)} \mathcal{v}_j^{(l)} (\mathbf{x}, {\bm{\theta}_j^{(l)}}).
\end{equation} 

We refer to \cref{eq:multilevel_fbpinn_network_architecture} as the multilevel FBPINN solution. Note, the original FBPINN solution described in \cref{sec:fbpinns} can be obtained by simply setting $L=1$; we refer to these as one-level FBPINNs going forward. 

Analogous to FBPINNs, we can train multilevel FBPINNs by using the same training scheme as PINNs and inserting \cref{eq:multilevel_fbpinn_network_architecture} into the PINN loss function. When using the hard-constrained PINN loss function~\cref{eq:pinn_loss_hard}, this yields the corresponding multilevel FBPINN loss function
\begin{equation} \label{eq:multilevel_fpinn_loss_hard} 
   \mathcal{L}(\bm{\theta}) =\frac{1}{N} \sum_{i=1}^{N} (\mathcal{N}[\mathcal{C} \frac{1}{L} \sum_{l=1}^L \sum_{j=1}^{J^{(l)}} \omega_j^{(l)} \mathcal{v}_j^{(l)}(\mathbf{x}_{i},\bm{\theta}_j^{(l)})] - f(\mathbf{x}_{i}))^{2}.
\end{equation}

\subsubsection{Example of a multilevel FBPINN}

We now show a simple example of a multilevel FBPINN to aid understanding. In particular, we use a multilevel FBPINN to solve the Laplacian boundary value problem,
$$
\begin{aligned}
- \Delta u & = f \quad \text{in } \Omega = [0,1]^d, \\
u & = 0 \quad \text{on } \partial\Omega.
\end{aligned}
$$

First we consider the 1D case ($d=1$), and set $f=8$. Then the exact solution is given by $u(x)=4x(1-x)$.

We create an $L=3$ level FBPINN to solve this problem, with $J^{(1)}=1$, $J^{(2)}=2,$ and $J^{(3)}=4$. Each level uses a uniform DD given by
\begin{align} \label{eq:used_dd}
\Omega_j^{(l)}=
\begin{cases}
\left[0.5-\delta/2,0.5+\delta/2\right] &l=1,\\
\left[\frac{(j-1)-\delta/2}{J^{(l)}-1},\frac{(j-1)+\delta/2}{J^{(l)}-1}\right] &l>1,
\end{cases}
\end{align}
where $\delta$ is defined as the \textit{overlap ratio} and is fixed at a value of $\delta=1.9$. Note that an overlap ratio of less than 1 means that the subdomains are no longer overlapping. The subdomain window functions form a partition of unity for each level and are given by
$$
\omega_j^{(l)} = \frac{\hat \omega_j^{(l)}}{\sum_{j=1}^{J^{(l)}}\hat \omega_j^{(l)}}
\quad
\text{where}
\quad
\hat \omega_j^{(l)}(x) = 
\begin{cases}
1 &l=1\\
[1+\cos(\pi(x-\mu_j^{(l)})/\sigma_j^{(l)})]^2 &l>1,
\end{cases}
$$
where $\mu_j^{(l)}=(j-1)/(J^{(l)}-1)$ and $\sigma_j^{(l)}=(\delta/2)/(J^{(l)}-1)$ represent the center and half-width of each subdomain respectively. The window functions for each level are plotted in \cref{fig:dd-collocation_points}~(a). A FCN \cref{eq:fcn} with 1 hidden layer, 16 hidden units, and $\tanh$ activation functions is placed in each subdomain, and the $x$ inputs to each subdomain network are normalized to the range [-1,1] over their individual subdomains.

\cor{Note that multilevel FBPINNs are not restricted to the particular choice of DD, level structure, window function, partition of unity, and subdomain network architecture used above. This choice may not be optimal, and the optimal choice is clearly problem-dependent. For example, it may be beneficial to use an irregular DD, level structure, and varying subdomain network sizes for problems where the solution has varying complexity in different parts of the domain.}

The multilevel FBPINN is trained using the hard-constrained loss function \cref{eq:multilevel_fpinn_loss_hard} with a constraining operator given by $[\mathcal{C} u] (x,\bm{\theta}) =\tanh(x/\sigma)\tanh((1-x)/\sigma)u(x,\bm{\theta})$ and $\sigma=0.2$. The loss function is minimized using the Adam optimizer with a learning rate of $1\times10^{-3}$ and $N=80$ uniformly-spaced collocation points across the domain.

The resulting multilevel FBPINN solution is shown in \cref{fig:dd-collocation_points}~(c), and the individual subdomain network solutions (with the constraining operator and window function applied) are shown in \cref{fig:dd-collocation_points}~(b). In this case, we find the FBPINN closely matches the exact solution.

Next we consider the 2D case ($d=2$), and set $f(x_1,x_2)=32(x_1(1-x_1)+x_2(1-x_2))$. Then the exact solution is given by $u(x_1, x_2) = 16(x_1(1-x_1)x_2(1-x_2))$.

In this case we create a $L=3$ level FBPINN to solve this problem, using a uniform rectangular DD for each level with $J^{(1)}=1\times1=1$, $J^{(2)}=2 \times 2=4$, and $J^{(3)}=4 \times 4=16$, as shown in \cref{fig:dd-collocation_points}~(d) and (e). The size of each subdomain along each dimension is defined similar as in~\cref{eq:used_dd} using, again, an overlap ratio of $\delta=1.9$. The subdomain window functions are given by
\begin{equation} \label{eq:used_window_function_nd}
\omega_j^{(l)} = \frac{\hat \omega_j^{(l)}}{\sum_{j=1}^{J^{(l)}}\hat \omega_j^{(l)}},
\quad
\text{where}
\quad
\hat \omega_j^{(l)}(\mathbf{x}) = 
\begin{cases}
1 &l=1\\
\prod_i^d [1+\cos(\pi(x_i-\mu_{ij}^{(l)})/\sigma_{ij}^{(l)})]^2 &l>1,
\end{cases}
\end{equation}
where $\mu_{ij}^{(l)}$ and $\sigma_{ij}^{(l)}$ represent the center and half-width of each subdomain along each dimension, respectively. A FCN~\cref{eq:fcn} with 1 hidden layer, 16 hidden units, and $\tanh$ activation functions is placed in each subdomain, and the $\mathbf{x}$ inputs to each subdomain network are normalized to the range [-1,1] along each dimension over their individual subdomains.

Similar to above, the multilevel FBPINN is trained using the hard-constrained loss function \cref{eq:multilevel_fpinn_loss_hard}, using a constraining operator given by $$[\mathcal{C} u] (\mathbf{x},\bm{\theta}) =\tanh(x_1/\sigma)\tanh((1-x_1)/\sigma)\tanh(x_2/\sigma)\tanh((1-x_2)/\sigma)u(\mathbf{x},\bm{\theta}),$$ with $\sigma=0.2$. The loss function is minimized using the Adam optimizer with a learning rate of $1\times10^{-3}$ and $N=80\times80=$~6,400 uniformly-spaced collocation points across the domain.

The resulting multilevel FBPINN solution is shown in \cref{fig:dd-collocation_points}~(f). Similar to the 1D case, we find the multilevel FBPINN solution closely matches the exact solution.

\subsubsection{Multilevel FBPINNs versus classical multilevel DDMs} \label{sec:classical-comparison}

Whilst multilevel FBPINNs are inspired by classical multilevel DDMs, a number of differences and similarities exist between these approaches. We believe it is insightful to briefly discuss these below.

Most classical DDMs can be described in terms of the abstract Schwarz framework~\cite{smith_domain_1996,toselli_domain_2005}. Similar to FBPINNs, this framework is based on a decomposition of a global function space $V$ into local spaces $\left\lbrace V_j \right\rbrace_{j=1}^J$ defined on overlapping subdomains $\Omega_j$, where
\begin{equation} \label{eq:schwarz:space_dec}
V = \sum_{j=1}^J R_j^\top V_j.
\end{equation}
Here, 
$R_j^\top: V_j \to V$ is an interpolation respectively prolongation operator from the local into the global space. These notions can be defined in a similar fashion at the continuous and discrete level. For the sake of simplicity, we suppose here a variational discretization of the PDE to solve. The space decomposition~\cref{eq:schwarz:space_dec} allows for decomposing any given discrete function $u \in V$ as
\begin{equation} \label{eq:schwarz:dec}
u = \sum_{j=1}^J R_j^\top v_j, \quad v_j \in V_j;
\end{equation}
due to the overlap, this decomposition is generally not unique. Schwarz DDMs are then based on solving local overlapping problems corresponding to the local spaces $\left\lbrace V_j \right\rbrace_{j=1}^J$ and merging them via the prolongation operators $R_j^\top$. 


Classical one-level Schwarz methods based on this framework are typically not scalable to large numbers of subdomains. In particular, since information is only transported via the overlap, their rate of convergence will deteriorate when increasing the number of subdomains~\cite{toselli_domain_2005}. In order to fix this, multilevel methods add coarser problems to the Schwarz framework to facilitate the global transfer of information; in particular, the coarsest level typically corresponds to a global problem.


We note that:
\begin{itemize}
\item In classical Schwarz methods, the global discretization space $V$ is often fixed first, and then, the local spaces $\left\lbrace V_j \right\rbrace_{j=1}^J$ are constructed. In FBPINNs, we do the opposite; we define a local space of neural network functions on an overlapping DD $\left\lbrace \Omega_j \right\rbrace_{j=1}^J$ and construct the global discretization space from them.

\item In classical Schwarz methods, the local functions $v_j \in V_j$ are generally not defined on the global domain $\Omega$ outside the overlapping subdomain $\Omega_j$; the prolongation operators $R_j^\top$ extend the local functions to $\Omega$ such that ${\rm supp}(R_j^\top v_j) \subset \Omega_j, \forall v_j \in V_j$. On the other hand, in FBPINNs, the local neural network functions $\mathcal{v}_j$ generally have global support, and the window functions $\omega_j$ are used to confine them to their subdomains. This difference \corrr{stems} from the fact that the local neural networks are not based on a spatial discretization but a function approximation; cf.~\cref{sec:pinns}. Nonetheless, both the prolongation operators and the window functions ensure locality; cf.~\cref{eq:schwarz:dec,eq:fbpinn_network_architecture}. Note that the prolongation operators in the restricted additive Schwarz (RAS) method~\cite{Cai:1999:RAS} also include a partition of unity, such that they are very close to the window functions in FBPINNs.

\item  A key difference is how the boundary value problem is solved. Whereas in DDMs, local subdomain problems are explicitly defined and solved in a global iteration, in FBPINNs, the global loss function is minimized. Moreover, classical DDMs can exploit properties of the system to be solved. For instance, if the PDE is linear elliptic, convergence guarantees for classical DDMs can be derived; cf.~\cite{toselli_domain_2005,Dolean:AID:2015}. In FBPINNs, we always have to solve a non-convex optimization problem~\cref{eq:pinn_loss_soft} or~\cref{eq:pinn_loss_hard}, which makes the derivation of convergence bounds difficult. Note that there are also nonlinear overlapping DDMs, for instance, additive Schwarz preconditioned inexact Newton (ASPIN) ~\cite{cai_nonlinearly_2002} and additive Schwarz preconditioned exact Newton (RASPEN) methods ~\cite{dolean:2016:NLP}.
\end{itemize}

\section{Numerical results} \label{sec:results}


In this section we assess the performance of multilevel FBPINNs. In particular, we investigate the accuracy and computational cost of using multilevel FBPINNs to solve various differential equations, and compare them to PINNs, PINNs with Fourier input features, \cor{SA-PINNs}, and one-level FBPINNs. 

First, in~\cref{sec:problems}, we introduce the problems studied. Then, in~\cref{sec:scaling} we introduce a notion of strong and weak scaling, inspired by classical DDMs, for assessing how the accuracy of FBPINNs and PINNs scales with computational effort and solution complexity. In~\cref{sec:common_implementation}, we list the common implementation details used across all experiments. Finally, in~\cref{sec:results_problems} we present our numerical results.

\subsection{Problems studied} \label{sec:problems}

The following problems are used to assess the performance of multilevel FBPINNs;

\subsubsection{Homogeneous Laplacian problem in two dimensions}

First, we consider the 2D homogeneous Laplacian problem already presented above, namely
\begin{equation} \label{eq:laplace_problem}
\begin{split}
- \Delta u & = f \quad \text{in } \Omega = [0,1]^2, \\
u & = 0 \quad \text{on } \partial\Omega,
\end{split}
\end{equation}
where
$$f(x_1,x_2)=32(x_1(1-x_1)+x_2(1-x_2)).$$
In this case, the exact solution is given by 
$$u(x_1, x_2) = 16(x_1(1-x_1)x_2(1-x_2)).$$

This problem is used to carry out simple ablation tests of the multilevel FBPINN. In particular, we assess how varying the number of levels and subdomains as well as the overlap ratio and size of the subdomain networks (architecture) affects the multilevel FBPINN performance.

\subsubsection{Multi-scale Laplacian problem in two dimensions}

Next, we consider a multi-scale variant of the Laplacian problem \cref{eq:laplace_problem} above by using the source term
\begin{equation} \label{eq:multiscale_laplace_f}
f(x_1,x_2) = \frac{2}{n}\sum_{i=1}^n \big( \omega_i \pi \big)^2 \sin({\omega_i \pi x_1}) \, \sin({\omega_i \pi x_2}).
\end{equation}
Then, the exact solution is given by 
$$
u(x_1,x_2) = \frac{1}{n}\sum_{i=1}^n \sin({\omega_i \pi x_1}) \, \sin({\omega_i \pi x_2}).
$$

In this case, multi-scale frequencies are contained in the solution, and the values of $n$ and $\omega_i$ allow us to control the number of components and the frequency of each component. We use this problem to assess how the performance of the multilevel FBPINN scales when more multi-scale components are added to the solution.

\subsubsection{Helmholtz problem in two dimensions}

Lastly, we study the 2D Helmholtz problem
\begin{equation} \label{eq:helmholtz_problem}
\begin{split}
\Delta u - k^2 u & = f \quad \text{in } \Omega = [0,1]^2, \\
u & = 0 \quad \text{on } \partial\Omega, \\
f(\mathbf{x}) & = e^{-\frac{1}{2}(\|\mathbf{x}-0.5 \|/\sigma)^2},
\end{split}
\end{equation}
with a constant (scalar) wave number, $k$. Here, homogeneous Dirichlet boundary conditions and a Gaussian point source with a scalar width, $\sigma$, placed in the center of the domain are used. Note that, for this problem, the exact solution is not known, and instead, we compare our models to the solution obtained from FD modeling, as described in~\cref{sec:helmholtz_fd}.

In this case, the solution contains complex patterns of standing waves where the dominant frequency of the solution depends on the wave number, $k$. We use this problem to test the multilevel FBPINN on a more realistic problem. We first carry out some simple ablation tests by assessing how varying the number of levels, subdomains, overlap ratio and size of the subdomain networks affects the multilevel FBPINN performance. Then, we assess how the performance of the multilevel FBPINN scales when the value of $k$ is increased.

\subsection{Definition of strong and weak scaling} \label{sec:scaling}


For both the multi-scale Laplacian and Helmholtz problems, we carry out strong and weak scaling tests. These assess how the accuracy of the multilevel FBPINN scales with computational effort and solution complexity and are inspired by the strong and weak scaling tests commonly used in classical DD. They are defined in the following way;
\begin{itemize}
\item \textit{Strong scaling:} We fix the complexity of the problem and increase the model capacity. For optimal scaling, we expect the convergence rate and/or accuracy to improve at the same rate as the increase of model capacity.
\item \textit{Weak scaling:} We increase the complexity of the problem and the model capacity at the same rate. For optimal scaling, we expect the convergence rate and/or accuracy to stay approximately constant.
\end{itemize}

For all our tests, increasing the model capacity means increasing the number of levels, number of subdomains, and/or the size of the subdomain networks. The exact factors varied and their rates of increase are detailed in the relevant results sections below. Note all of the multilevel FBPINNs tested have been trained on a single GPU, and hence we only show strong and weak scaling tests with respect to model capacity and not hardware parallelization.

\subsection{Common implementation details} \label{sec:common_implementation}

\begin{figure}[!t]
\centering
\includegraphics[width=0.75\textwidth]{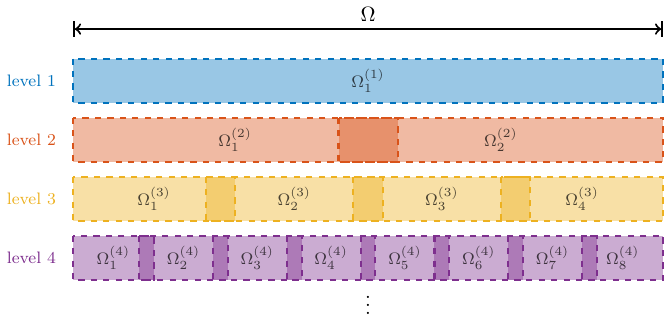}
\caption{Hierarchy of levels used in the multilevel FBPINN. For all the multilevel FBPINNs tested we use an exponential level structure. This means that the number of subdomains in each level is given by $2^{d(l-1)}$, where $l$ is the level number and $d$ is the dimensionality of the domain. Our hypothesis is that this helps the multilevel FBPINN model solutions with frequency components that span multiple orders of magnitude. 
\label{fig:multilevel-fbpinns}}
\end{figure}


Many of the implementation details of the multilevel FBPINNs, one-level FBPINNs, and PINNs tested are the same across all tests. These details are presented here; some are only changed for ablation studies, in which case they are described in the relevant results section below.

\paragraph{Level structure} Firstly, all multilevel FBPINNs use an exponentially increasing number of subdomains per level. In particular, we choose $J^{(l)}=2^{d(l-1)}$ for $l=1, ..., L$. This level structure is shown in \cref{fig:multilevel-fbpinns}. This constraint is chosen so that the multilevel FBPINN is able to contain an exponentially large number of subdomains with a relatively small number of levels; our hypothesis is that this helps the multilevel FBPINN model solutions with frequency components that span multiple orders of magnitude. 

\paragraph{Domain decomposition} All FBPINNs tested use a uniform rectangular DD for each level, with all multilevel FBPINNs having $2^{l-1}$ subdomains along each dimension. The size of each subdomain along each dimension is defined similar to \cref{eq:used_dd}, i.e., all 2D DDs look similar to those shown in \cref{fig:dd-collocation_points}~(d) and (e). Furthermore, all FBPINNs use the same subdomain window functions, given by \cref{eq:used_window_function_nd}.

\paragraph{Network architecture} All FBPINNs tested use FCNs with identical architectures as their subdomain networks. The PINNs tested either use FCNs or FCNs with Fourier input features (see \cref{sec:fourierfeatures} for the definition of Fourier features). For all the FBPINNs tested, the $\mathbf{x}$ inputs to each subdomain network are normalized to the range [-1,1] along each dimension over their individual subdomains. For the PINNs tested, the $\mathbf{x}$ inputs are normalized to the range [-1,1] along each dimension over the global domain. $\tanh$ is used for all activation functions.

\paragraph{Loss function and optimization} All FBPINNs and PINNs tested use the hard-constrained variants of their loss functions. For fairness, the same constraining operator is used across all models tested for a given problem.  Furthermore, the same collocation points are used for training whenever multiple models are compared on a given problem. This is similarly the case for all testing points used after training. \cor{The SA-PINNs tested use learnable weights for each collocation point in their loss function; this is described in more detail in \cref{sec:sapinns}}. All tests use the Adam optimizer with a learning rate of $1\times10^{-3}$ \corr{ except for the PINNs with Fourier input features, which are trained using a learning rate of $1\times10^{-4}$, because it was found their convergence is unstable when using a larger learning rate}. For robustness, all models are trained 10 times using different random starting seeds, and all results are reported as averages over these different seeds. All models are evaluated using the normalized L1 test loss, given by $\mathcal{L}(\bm{\theta})=\frac{1}{M}\sum_i^M \| u(\mathbf{x}_i,\bm{\theta})-u(\mathbf{x}_i) \| / \sigma$, where $M$ is the number of test points and $\sigma$ is the standard deviation of the set of true solutions $\{u(\mathbf{x}_i)\}^M_i$.

\paragraph{Software and hardware implementation} All FBPINNs and PINNs tested are implemented using a common training framework written in JAX \cite{jax}. Further details on our software implementation are given in \cref{sec:software}. All models are trained on a single NVIDIA RTX 3090 GPU.

\subsection{Results} \label{sec:results_problems}

Here, we will discuss the results for the model problems described in~\cref{sec:problems}.

\begin{figure}[!t]
\centering
\includegraphics[width=0.98\textwidth]{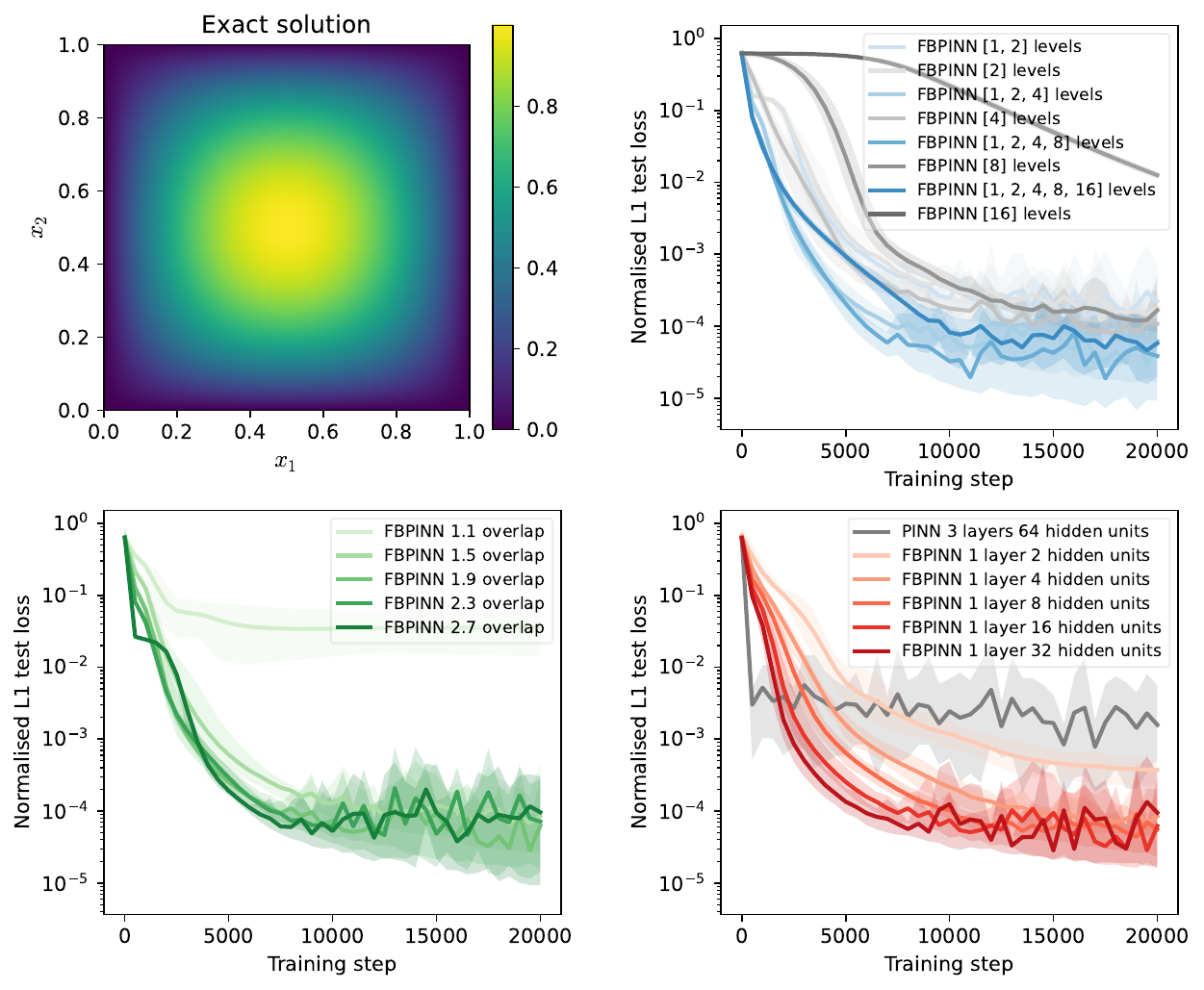}
\caption{Ablation tests using the homogeneous Laplacian problem. The convergence curve of a baseline multilevel FBPINN is plotted when changing the number of levels (top right), overlap ratio (bottom left), and number of hidden units for each subdomain network (bottom right). The baseline model has $L=3$ levels, an overlap ratio of $\delta=1.9$, and 16 hidden units for each subdomain network. The exact solution is shown (top left). Convergence curves of two other benchmarks are shown; a PINN (bottom right), and one-level FBPINNs with varying numbers of subdomains (top right). The lists which label each model in the top right plot contain the number of subdomains along each dimension for each level in the model. Filled region edges show the minimum and maximum loss values across 10 random starting seeds and lines show the average.
\label{fig:ablation_laplace}}
\end{figure}

\subsubsection{Homogeneous Laplacian problem in two dimensions} \label{sec:homogeneous_laplace}

First, we carry out simple ablation tests of the multilevel FBPINN using the 2D homogeneous Laplacian problem described in \cref{sec:homogeneous_laplace}. 

To carry out our ablation tests, we first train a baseline multilevel FBPINN to solve this problem, using $L=3$ levels, an overlap ratio along each dimension of $\delta=1.9$, and FCNs with 1 hidden layer and 16 hidden units for each subdomain network. The multilevel FBPINN is trained using the constraining operator 
$$
    [\mathcal{C} u] (\mathbf{x},\bm{\theta}) =\tanh(x_1/\sigma)\tanh((1-x_1)/\sigma)\tanh(x_2/\sigma)\tanh((1-x_2)/\sigma)u(\mathbf{x},\bm{\theta})
$$ 
with $\sigma=0.2$. Here, $N=80\times80=$~6,400 uniformly-spaced collocation points and $M=350\times 350$ uniformly-spaced test points across the global domain are used to train and test the multilevel FBPINN, respectively.

Given this baseline model, we then vary different hyperparameters over a range of values and measure the change in performance. This is carried out for the number of levels ranging from $L=2$ to $5$, the overlap ratio ranging from $1.1$ to $2.7$, and the number of hidden units in the subdomain network ranging from $2$ to $32$. Our results are shown in~\cref{fig:ablation_laplace}. We observe that the accuracy of the multilevel FBPINN does not depend significantly on the number of levels, likely because in this case the solution is very simple. However, its accuracy increases as the overlap ratio increases, likely because there is more communication between the subdomain networks, which is similar to what is expected in classical DDMs. Furthermore, its accuracy increases as the number of free parameters of the subdomain networks increases. This is expected as the capacity of the model increases. Thus, the multilevel FBPINN has similar characteristics to classical DDMs for this problem.

We carry out two other benchmark tests. First, we train a PINN with 3 hidden layers and 64 hidden units, and second, we train four one-level FBPINNs with $J^{(1)}=2,4,8, \text{ and } 16$ subdomains along each dimension, respectively. All other relevant hyperparameters are kept the same as the baseline model. These results are also shown in~\cref{fig:ablation_laplace}. In these tests, the PINN is able to solve the problem, although its final accuracy is lower than the baseline multilevel FBPINN and its convergence curve is more unstable. Furthermore, the accuracy of the one-level FBPINN reduces as more subdomains are added. This is analogous to the expected behavior of one-level classical DDMs, which is not scalable to large numbers of subdomains, and shows that coarse levels are required for scalability. It is therefore likely that the additional levels in FBPINNs serve the same purpose as in classical DDMs, i.e., they allow direct transfer of global information.

\subsubsection{Multi-scale Laplacian problem in two dimensions} \label{sec:multi-scale_laplace}

\begin{figure}[!t]
\centering
\includegraphics[width=\textwidth]{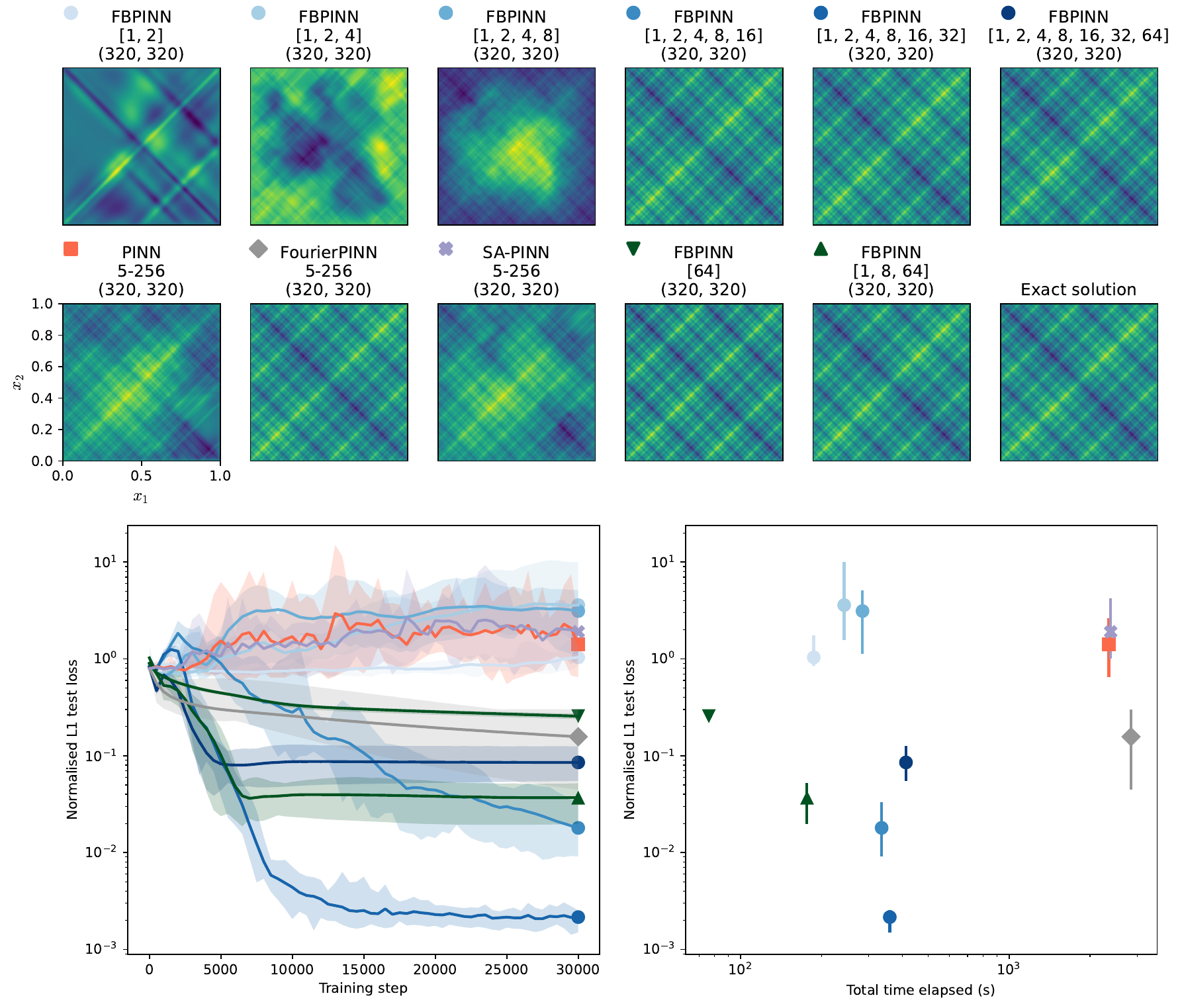}
\caption{Strong scaling test using the multi-scale Laplacian problem. In this test the problem complexity is fixed and the solution estimated using multilevel FBPINNs with increasing numbers of levels are plotted (top row). The title of each plot describes the level structure (first line) and the number of collocation points along each dimension (second line). The color-coded convergence curves and training times for each model are shown (bottom row). Filled region edges show the minimum and maximum loss values across 10 random starting seeds and lines show the average. Error bars show the minimum and maximum loss values and training times. The exact solution is shown (middle row). Plots of the solutions and convergence curves of a PINN, PINN with Fourier input features, \cor{SA-PINN}, one-level FBPINN and three-level FBPINN benchmark are also shown (middle and bottom row).
\label{fig:strong_laplace}}	
\end{figure}

\begin{figure}[!t]
\centering
\includegraphics[width=\textwidth]{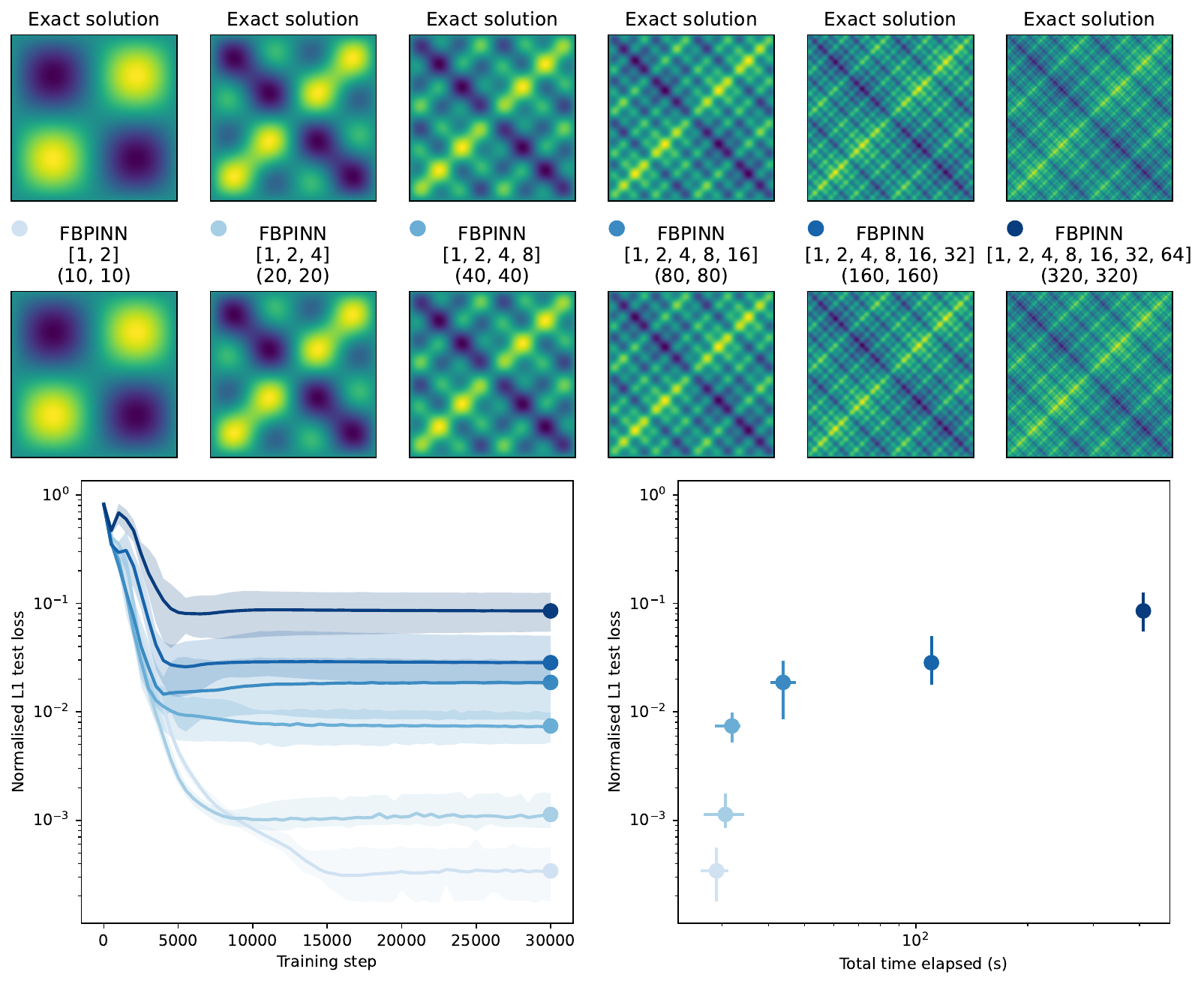}
\caption{Weak scaling test using the multi-scale Laplacian problem. In this test the problem complexity is increased (in this case, the number of frequency components in the solution) (top row) and the solution estimated using multilevel FBPINNs with increasing numbers of levels and collocation points are plotted (middle row). The title of each plot describes the level structure (first line) and the number of collocation points along each dimension (second line). The color-coded convergence curves and training times for each model are shown (bottom row). Filled region edges show the minimum and maximum loss values across 10 random starting seeds and lines show the average. Error bars show the minimum and maximum loss values and training times. 
\label{fig:weak_laplace}}
\end{figure}

Next, we evaluate the strong and weak scalability of the multilevel FBPINN using the multi-scale Laplacian problem described in~\cref{sec:multi-scale_laplace}. 

\paragraph{Strong scaling test} First, we carry out a strong scaling test. Here, the problem complexity is fixed and we assess how the performance of the multilevel FBPINN changes as the capacity of the model is increased. In particular, we fix the problem complexity by choosing $n=6$ with $\omega_i=2^i$ for $i=1, ..., n$ in \cref{eq:multiscale_laplace_f}. Thus, the solution contains 6 multi-scale components with exponentially increasing frequencies. This represents a much more challenging problem than the homogeneous problem studied above. The exact solution in this case is shown in \cref{fig:strong_laplace}.

We then increase the capacity of the multilevel FBPINN by increasing the number of levels, testing from $L=2$ to $7$. The rest of the hyperparameters of the multilevel FBPINN are kept fixed across all tests. Namely, we use $N=320\times320=$~102,400 uniformly-spaced collocation points throughout the domain, an overlap ratio of $\delta=1.9$ and FCNs with 1 hidden layer and 16 hidden units for each subdomain network. All models are trained using the constraining operator $[\mathcal{C} u] (\mathbf{x},\bm{\theta}) =\tanh(x_1/\sigma)\tanh((1-x_1)/\sigma)\tanh(x_2/\sigma)\tanh((1-x_2)/\sigma)u(\mathbf{x},\bm{\theta})$ with $\sigma=1/\omega_n$. $M=350\times 350$ uniformly-spaced test points are used to test all models.

The results of this study are shown in~\cref{fig:strong_laplace}. We find that the accuracy of the multilevel FBPINN increases as the number of levels increases, where the $L=2,3, \text{ and } 4$ models are unable to accurately model the solution, whilst the $L=5,6 \text{ and } 7$ models are able to accurately model all of the frequency components. The test shows that the multilevel FBPINN is able to solve a high frequency, multi-scale problem, and exhibits strong scaling behavior somewhat analogous to what is expected by classical DDMs. \corrr{However, we note that the accuracy of the 7-level FBPINN is worse than the 6-level FBPINN. We believe this may because at the finest level of the 7-level PINN, each subdomain network only contains approximately $10\times10$ collocation points. More collocation points may allow this level to converge more accurately.}

\cor{Five} other benchmark tests are carried out for this problem. First, we train a PINN with 5 hidden layers and 256 hidden units, a PINN with 256 Fourier input features \corr{with $\sigma=5$}, 5 hidden layers and 256 hidden units, \cor{and a SA-PINN with 5 hidden layers and 256 hidden units}. Then, we train a one-level FBPINN with $J^{(1)}=64$ subdomains along each dimension and a three-level FBPINN with $J^{(1)}=1$, $J^{(2)}=8$, and $J^{(3)}=64$ subdomains along each dimension, respectively. All other relevant hyperparameters are kept the same as the baseline model above. These results are also shown in~\cref{fig:strong_laplace}. We find that the accuracy of the \cor{standard} PINN is poor, and it is only able to model some of the cycles in the solution. Furthermore its convergence curve is very unstable, and its training time is an order of magnitude larger than the $L=7$ level FBPINN tested. Its poor convergence is likely due to spectral bias and the increasing complexity of the PINN's optimization problem, as discussed in \cref{sec:fbpinns} and \cite{Moseley2023}. For this problem, adding Fourier features to the PINN significantly improves its accuracy, although it its training time remains high and it converges slower than the multilevel FBPINNs tested. \cor{The SA-PINN does not offer any improvement over the standard PINN.} The one-level FBPINN is able to model the solution, although its accuracy is less than the $L=7$ level FBPINN. Finally, the three-level FBPINN benchmark performs similarly to the $L=7$ level FBPINN. This suggests that multilevel FBPINNs with stronger coarsening ratios can be used, and that multilevel FBPINNs are not strongly dependent on their coarsening ratio.

\paragraph{Weak scaling test} Next, we carry out a weak scaling test. Here, the problem complexity, number of collocation points and model capacity are scaled at the same rate, and we assess how the performance of the multilevel FBPINN changes. We increase the model capacity in the same way as the strong scaling test above, i.e., the number of levels is increased from $L=2$ to $7$. However, now the problem complexity is also scaled, such that for each test $n=L-1$ and $\omega_i=2^i$ for $i=1, ..., n$. Furthermore, each test has $(5 \times 2^{L-1}) \times (5 \times 2^{L-1})$ uniformly-spaced collocation points. Note that the number of subdomains, number of collocation points, and the frequency range of the solution all grow exponentially, and the multilevel FBPINN is in alignment with the problem structure. All other hyperparameters are fixed to the same values as the strong scaling test above.

The results of this test are shown in \cref{fig:weak_laplace}. We find that the multilevel FBPINNs are able to model all of the problems tested accurately, that is, modeling all of their frequency components. However, the normalized L1 accuracy of the multilevel FBPINNs does reduce somewhat as the problem complexity increases. Thus in this test the multilevel FBPINN exhibits near -- but not perfect -- weak scaling. \corrr{As an additional test we plot the contribution to the FBPINN solution from each level for the 5-level FBPINN in \cref{sec:levels}.}

\subsubsection{Helmholtz problem in two dimensions} \label{sec:helmholtz}

\begin{figure}[!t]
\centering
\includegraphics[width=0.98\textwidth]{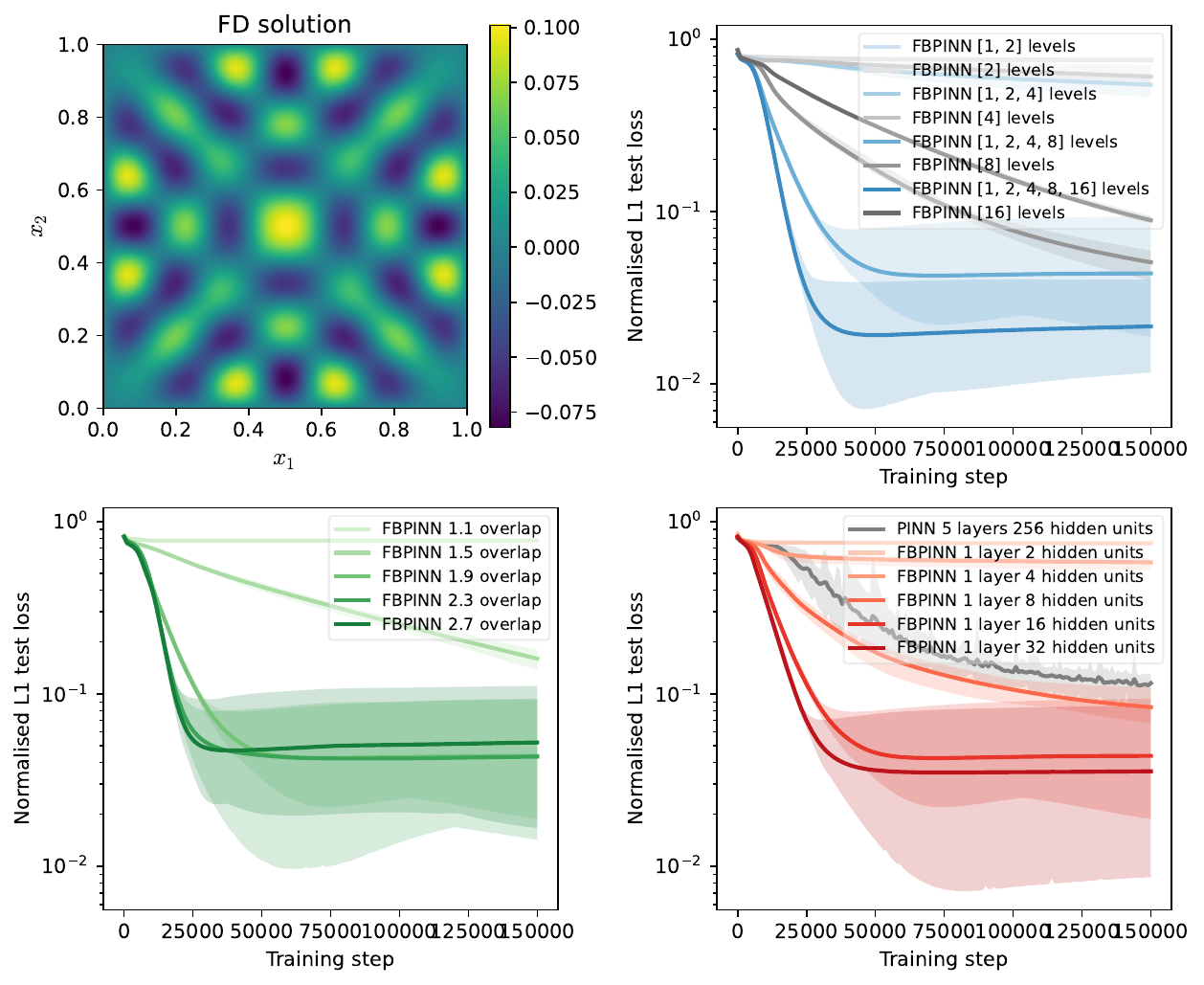}
\caption{Ablation tests using the Helmholtz problem. The convergence curve of a baseline multilevel FBPINN is plotted when changing the number of levels (top right), overlap ratio (bottom left), and number of hidden units for each subdomain network (bottom right). The baseline model has $L=4$ levels, an overlap ratio of $\delta=1.9$, and 16 hidden units for each subdomain network. The solution obtained from FD modeling is shown (top left). Convergence curves of two other benchmarks are shown; a PINN (bottom right), and one-level FBPINNs with varying numbers of subdomains (top right). The lists which label each model in the top right plot contain the number of subdomains along each dimension for each level in the model. Filled region edges show the minimum and maximum loss values across 10 random starting seeds and lines show the average.
\label{fig:ablation_helmholtz}}
\end{figure}

\begin{figure}[!t]
\centering
\includegraphics[width=0.94\textwidth]{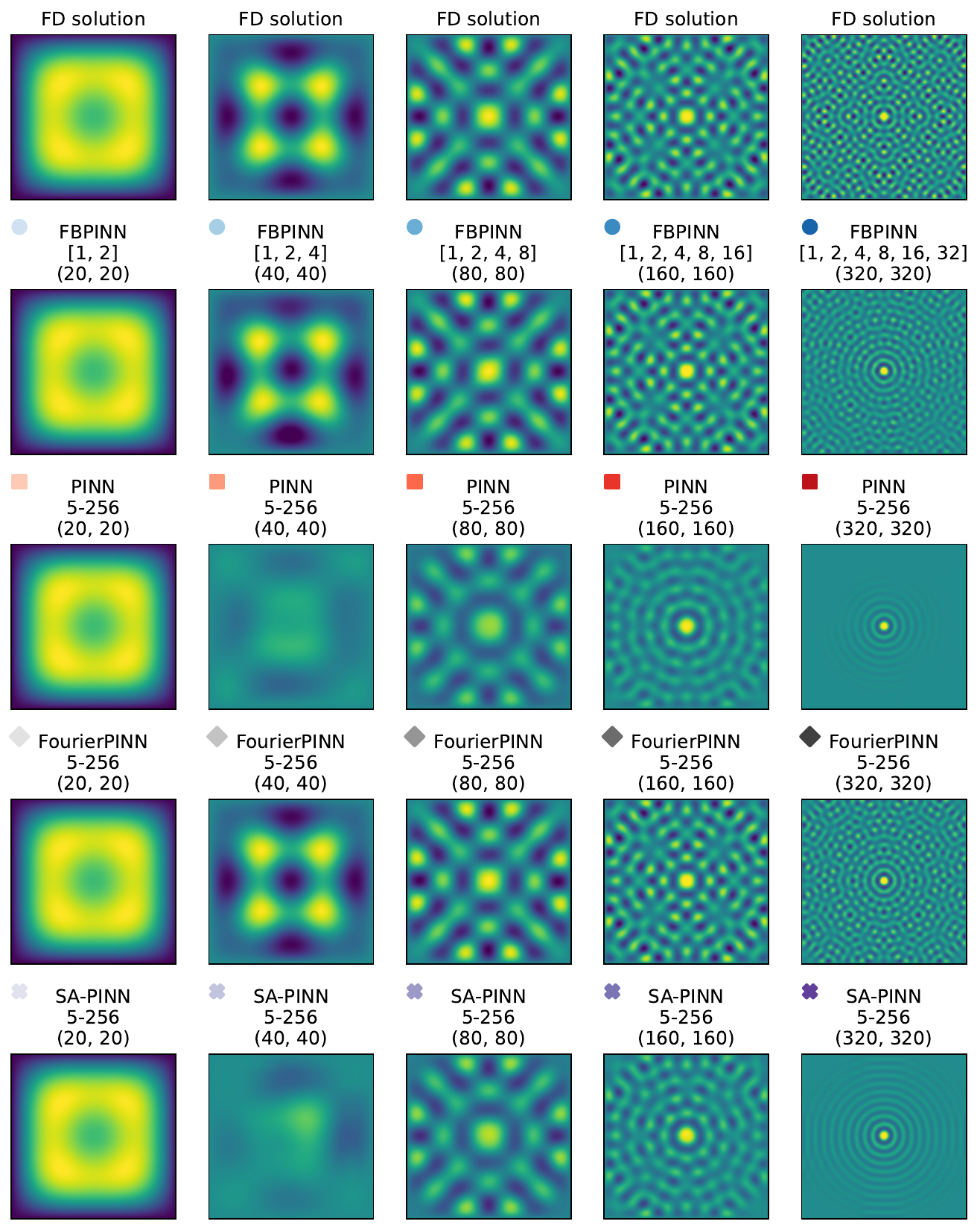}
\caption{Weak scaling test using the Helmholtz problem. In this test the problem complexity is increased (in this case, the wave number) (top row) and the solution estimated using multilevel FBPINNs with increasing numbers of levels and collocation points are plotted (second row). The title of each plot describes the level structure (first line) and the number of collocation points along each dimension (second line). \cor{Three} benchmarks using a PINN, a PINN with Fourier input features, \cor{and a SA-PINN, all} with a fixed network size and increasing numbers of collocation points are also shown (third and fourth row).
\label{fig:weak_helmholtz_constant_solution}}
\end{figure}
\clearpage

\begin{figure}[!t]
\centering
\includegraphics[width=1\textwidth]{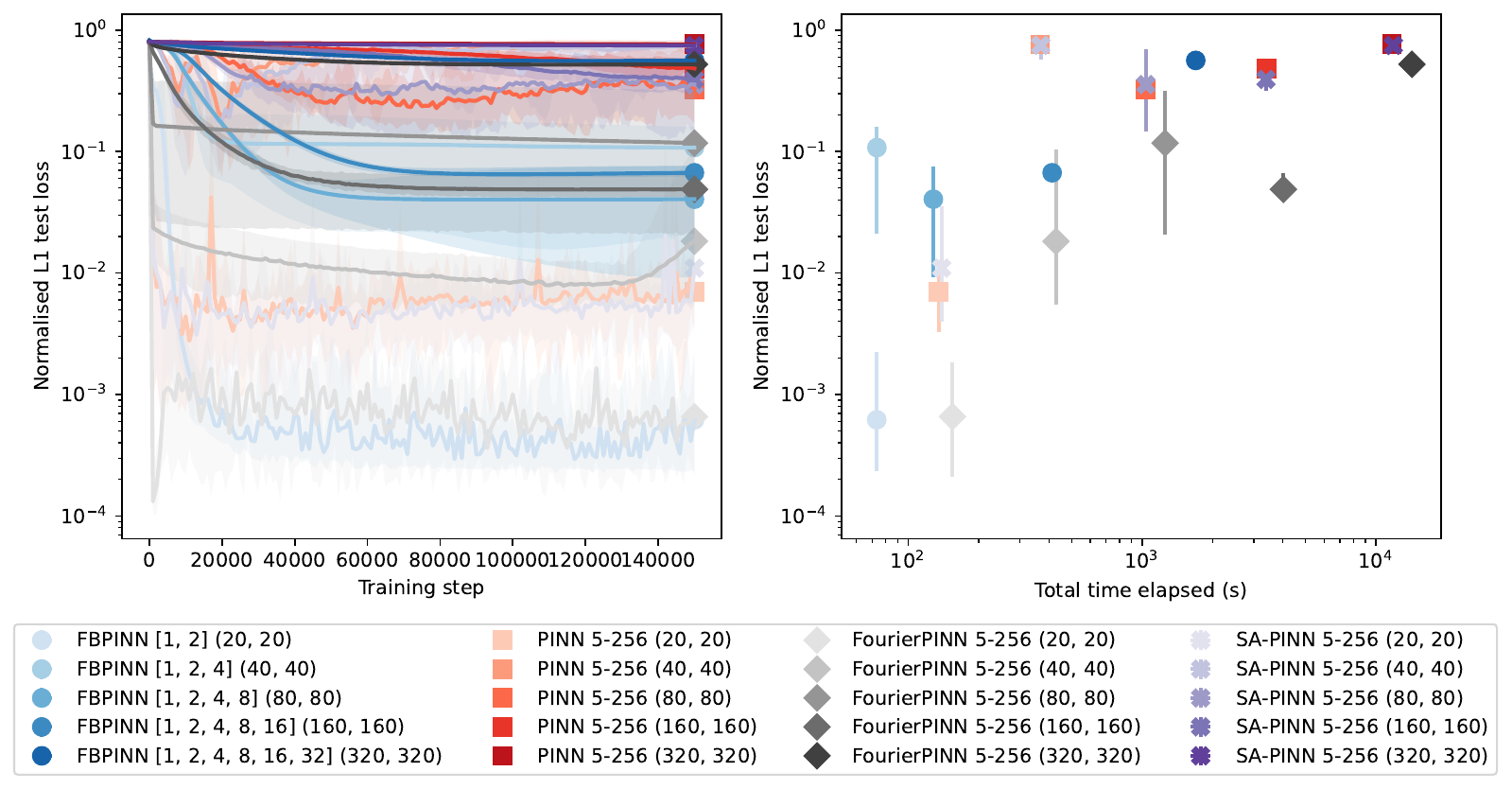}
\caption{Color-coded convergence curves and training times for each model displayed in \cref{fig:weak_helmholtz_constant_solution}. Filled region edges show the minimum and maximum loss values across 10 random starting seeds and lines show the average. Error bars show the minimum and maximum loss values and training times.
\label{fig:weak_helmholtz_constant_loss}}
\end{figure}

Finally, we test the multilevel FBPINN using the more complex Helmholtz problem described in \cref{sec:helmholtz}. Again, we carry out ablation tests first and then carry out a weak scaling study assessing how the performance of the multilevel FBPINN changes as the wave number, $k$, increases.

\paragraph{Ablation tests} For our ablation tests, we fix the problem parameters to be $k=2^4 \pi / 1.6$ and $\sigma=0.8/2^4$ in \cref{eq:helmholtz_problem}. Then, similar to \cref{sec:homogeneous_laplace}, we train a baseline multilevel FBPINN to solve this problem, using $L=4$ levels, an overlap ratio along each dimension of $\delta=1.9$, and FCNs with 1 hidden layer and 16 hidden units for each subdomain network. The multilevel FBPINN is trained using the constraining operator $[\mathcal{C} u] (\mathbf{x},\bm{\theta}) =\tanh(x_1/\sigma)\tanh((1-x_1)/\sigma)\tanh(x_2/\sigma)\tanh((1-x_2)/\sigma)u(\mathbf{x},\bm{\theta})$ with $\sigma=1/k$. We use $N=160 \times 160=$~25,600 uniformly-spaced collocation points and $M=320\times 320$ uniformly-spaced test points to train and test the multilevel FBPINN, respectively.

Given this baseline model, we then vary different hyperparameters over a range of values and measure the change in performance. This is carried out for the number of levels ranging from $L=2$ to $5$, the overlap ratio ranging from $1.1$ to $2.7$, and the number of hidden units in the subdomain network ranging from $2$ to $32$. Our results are shown in~\cref{fig:ablation_helmholtz}. We obtain similar results to the ablation tests carried out in~\cref{sec:homogeneous_laplace} for the homogeneous Laplace problem. Namely, that the accuracy of the multilevel FBPINN improves as the overlap ratio and the number of free parameters of the subdomain networks increases. Furthermore, its accuracy improves as the number of levels increases, likely because the solution contains relatively high frequencies and multiple subdomains are needed. We observe relatively large variations of the test loss between different random initialisations of the models below a value of $10^{-1}$. In particular, the final loss can be somewhere between $10^{-2}$  and $10^{-1}$.

We carry out two other benchmark tests. First, we train a PINN with 5 hidden layers and 256 hidden units, and second, we train four one-level FBPINNs with $J^{(1)}=2,4,8, \text{ and } 16$ subdomains along each dimension, respectively. All other relevant hyperparameters are kept the same as the baseline model. These results are also shown in \cref{fig:ablation_helmholtz}. Here, the PINN converges poorly, which again highlights the shortcomings of PINNs when solving more complex problems. Furthermore, the convergence of all the one-level FBPINNs is much slower than the multilevel FBPINN, and their final accuracy is worse. This again suggests that multiple levels are required for scalability.

\paragraph{Weak scaling test} We carry out a weak scaling study, where both the problem complexity and model capacity are scaled at the same rate. In a similar fashion to the weak scaling test in~\cref{sec:multi-scale_laplace}, the capacity of the multilevel FBPINN is increased by increasing the number of levels, testing from $L=2$ to $6$. For each test, $(10 \times 2^{L-1}) \times (10 \times 2^{L-1})$ uniformly-spaced collocation points are used. The problem complexity for each test is increased by setting $k=2^L \pi / 1.6$ and $\sigma=0.8/2^L$ in \cref{eq:helmholtz_problem}. All other hyperparameters are fixed to the same values as the baseline model used in the ablation tests above.

The results of this test are shown in \cref{fig:weak_helmholtz_constant_solution} and \cref{fig:weak_helmholtz_constant_loss}. We find that the multilevel FBPINN is able to accurately model all the problems tested, except for the highest wave number test. In this case, the multilevel FBPINN successfully models the dominant frequency and overall concentricity of the solution but fails to model its more complex motifs. In this case, we believe that the FBPINN is struggling to satisfy both the point source and Dirichlet boundary conditions. Without the Dirichlet boundary condition, the solution to \cref{eq:helmholtz_problem} is that of a simple point source. For all tests, we notice that in the first few training steps this is the solution learned by the multilevel FBPINN, which is then updated to the correct solution after further training. Thus, it appears the presence of the Dirichlet boundary condition leads to an optimization problem which remains challenging. This is consistent with challenges arising in solving Helmholtz problem using classical iterative numerical solvers. Further work is required to understand this behaviour; one may be able to address this problem by using subdomain scheduling strategies to incrementally train the multilevel FBPINN, as proposed in~\cite{Moseley2023}.

Finally, we carry out the same weak scaling test but using a PINN instead of a multilevel FBPINN. For each test, the PINN's architecture is kept fixed at 5 hidden layers and 256 hidden units whilst the number of collocation points and problem complexity is increased in the same way as the previous test. We also test a PINN with 256 Fourier input features, 5 hidden layers and 256 hidden units, \cor{and a SA-PINN with 5 hidden layers and 256 hidden units}. \corr{The $\sigma$ values of the Fourier input features are hand-selected for each test, and are $0.4, 1, 2, 1.5,$ and $5$ in order of increasing problem complexity}. All other relevant hyperparameters are kept the same \corr{as the FBPINNs tested}. The result of this study is shown in~\cref{fig:weak_helmholtz_constant_solution} and \cref{fig:weak_helmholtz_constant_loss}. In this case, we find that \corr{the PINN and SA-PINN are unable to accurately model any of the solutions, and their training time is an order of magnitude larger than the multilevel FBPINN. The PINN with Fourier input features is able to model the solutions with a similar level of accuracy as the multilevel FBPINN, but its training time is an order of magnitude larger. Thus, multilevel FBPINNs still outperform PINNs for this problem.}

\section{Discussion} \label{sec:discussion}









Across all the problems studied, we find that the multilevel FBPINNs consistently outperform the one-level FBPINNs and PINNs tested. The multilevel FBPINNs are more accurate than the one-level FBPINNs when a large number of subdomains are used, suggesting that coarse levels are required for scalability by improving the global communication. Furthermore, the multilevel FBPINNs are significantly more accurate and computationally efficient than the PINNs tested. \cor{However, we have only started to investigate multilevel FBPINNs in this work and there are many important research questions outstanding.}

\cor{One important question would be to investigate the performance of multilevel FBPINNs on problems with complex geometries and solutions which have varying complexity in different parts of the domain. In this work, we restrict ourselves to problems with rectangular geometries and homogeneous solution complexity, and use multilevel FBPINNs with uniform rectangular decompositions and an exponential level structure. However, for problems with complex geometry and varying solution complexity, it is likely to be beneficial to use irregular DDs with irregular level structures and varying subdomain network sizes, to capture inhomogeneities in the solution. Whilst the FBPINN framework can represent such a model, implementing it is practically challenging for two reasons. First, to remain computationally efficient it is likely that a fully asynchronous training code across subdomains is required, because each subdomain network would require a different amount of computation to be evaluated and trained. Secondly, it may be challenging to choose an optimal DD, level structure and the subdomain network sizes, especially if characteristics of the solution are not known beforehand. One interesting direction here would be to try to learn the DD, for example by jointly learning the parameters of the FBPINN window functions with the subdomain network parameters or using a gating network similar to \cite{Stiller2020, Hu2022}.}


Another valuable direction would be to study the theoretical convergence properties of multilevel FBPINNs. A major limitation of PINNs compared to classical DDMs is that their convergence properties are still poorly understood. In particular, whilst the multilevel FBPINN exhibits good scaling properties for the Laplacian problems studied, it remains unclear why the optimization of the high wave number Helmholtz problem is challenging; note that the convergence of classical DDMs for high wave number Helmholtz problems is also not fully understood.


\cor{Furthermore, it is important to investigate ways to accelerate the computational efficiency of multilevel FBPINNs further. Whilst multilevel FBPINNs  are }over an order of magnitude more efficient than the PINNs tested, their training times are still likely to be slower than many traditional methods, such as numerical solvers for finite difference or finite element systems. Fundamentally, this is because (FB)PINNs yield a non-convex optimization problem, which is relatively expensive compared to the linear solves which traditional methods typically rely on. One way to accelerate (multilevel) FBPINNs, as suggested in \cite{Moseley2023}, is to provide more inputs to the subdomain networks, such as BCs and PDE coefficients, and train across a range of these inputs so that the multilevel FBPINN learns a fast surrogate model which does not need to be retrained for each new solution. 


\cor{Alongside this, multi-GPU training of (multilevel) FBPINNs should be investigated. Here we solve all problems using a single GPU, but multi-GPU training will become essential for problem sizes where $10,000+$ subdomains are required (for example, 3D problems, or problems with highly multi-scale solutions). In \cref{sec:fbpinn-efficiency}, we show that (multilevel) FBPINNs are theoretically scalable to large problem sizes: the computational cost of evaluating the FBPINN solution is $\mathcal{O}(NC \tilde S)$, where $C$ is the average number of subdomains a collocation point belongs to, $\tilde S$ is the cost of computing the output of a single subdomain network for a single collocation point and $N$ is the number of collocation points. Importantly, this cost is independent of the total number of subdomains ($J$), and scales linearly with the number of collocation points. Practically, it may be challenging to achieve perfect linear scaling when using multiple GPUs because of the communication required between GPUs. For FBPINNs, the only communication required between subdomains is within their overlapping regions, where subdomain solutions are summed together; note that, in the multilevel case, subdomains on different levels may also overlap, therefore, increasing the required communication but improving the numerical scalability. One possible parallel implementation of FBPINNs was proposed by \cite{Moseley2023} (Algorithm~1 and Figure~4), where a separate GPU is used to train each subdomain network. Importantly, communication between GPUs is only required in the forward pass when summing the outputs of the subdomain networks in the overlapping regions; the backpropagation and updating of the subdomain networks can then be done independently on each GPU for each subdomain network. Future work will investigate the parallel scalability of FBPINNs in detail.}

\section*{Data availability}
All the code for reproducing this work is available here: \href{https://github.com/benmoseley/FBPINNs/tree/multilevel-paper/multilevel-paper}{https://github.com/benmoseley/FBPINNs/tree/multilevel-paper/multilevel-paper}. All data is generated synthetically using this code.

\section*{CRediT contributions}
\noindent Conceptualization (VD, AH, SM, BM);
Data curation (BM);
Investigation (VD, AH, SM, BM);
Methodology (VD, AH, SM, BM);
Project administration (VD, AH, SM, BM);
Resources (BM, SM);
Software (BM);
Validation (BM);
Visualization (BM);
Writing - original draft (VD, AH, SM, BM);
Writing - review \& editing (VD, AH, SM, BM).


\appendix


\section{Software implementation} \label{sec:software}
%
%
All FBPINNs and PINNs are implemented using a common training framework written using the JAX automatic differentiation library~\cite{jax}. When training (multilevel) FBPINNs, computing the FBPINN solution (either~\cref{eq:fbpinn_network_architecture} or~\cref{eq:multilevel_fbpinn_network_architecture}) naively can be very expensive. This is because evaluating the solution at each collocation point involves summing over all subdomain networks and all levels. However, the cost of this summation can be significantly reduced by exploiting that, because the output of all subdomain networks is zero outside of the corresponding subdomains, only subdomains which contain each collocation point contribute to the summation. Practically, this can be carried out by pre-computing a mapping describing which subdomains contain each collocation point before training and only evaluating the corresponding subdomain networks during training. Another important efficiency gain in our software implementation is that the outputs of each subdomain network are computed in parallel on the GPU by using JAX's \verb|vmap| functionality. This is important as the FBPINNs tested use small subdomain networks that if evaluated sequentially would not fully utilize the GPU's parallelism.

\section{Finite difference solver for Helmholtz equation} \label{sec:helmholtz_fd}
We use a finite difference (FD) solver to compute a reference solution for the Helmholtz problem studied in \cref{sec:helmholtz}. For all the problem variants studied, we discretize the Laplacian operator in \cref{eq:helmholtz_problem} using a 5-point stencil, and we discretize the solution using a $320\times320$ uniformly-spaced mesh over the problem domain. This turns \cref{eq:helmholtz_problem} into a set of linear equations, which are solved using the \verb|scipy.sparse.linalg| \cite{scipy} sparse direct solver, that is, using UMFPACK~\cite{davis2004umfpack}.

\section{Definition of Fourier features}
\label{sec:fourierfeatures}

We compare the performance of FBPINNs to PINNs with Fourier features for a number of experiments. Fourier features can help neural networks learn high-frequency functions \cite{Tancik2020} and they have been shown to help the convergence of PINNs when solving multi-scale problems \cite{Wang2021d}. When using Fourier features, the inputs of neural networks are transformed using trigonometric functions before inputting them into the network. These Fourier features are given by
\begin{equation}
\gamma(x) = [\cos(2 \pi \Gamma x),\,\sin(2 \pi \Gamma x)]~,
\end{equation}
where $\Gamma$ is a matrix of shape $k\times d$, $k$ is the number of Fourier features and $x$ is the input vector $x \in \mathbb{R}^{d}$. The values of $\Gamma$ represent the frequency of the features and they are typically sampled from a univariate Gaussian distribution with a mean and standard deviation denoted by $\mu$ and $\sigma$. \corr{Similar to \cite{Tancik2020, Wang2021d}, we fix $\mu=0$ and hand-select $\sigma$ based on test accuracy for all tests in this work. We find that the convergence of the PINN is highly sensitive to the value of $\sigma$ and the learning rate.}

\begin{core}

\section{Definition of self-adaptive PINNs (SA-PINNs)}
\label{sec:sapinns}

We compare the performance of FBPINNs to self-adaptive SA-PINNs~\cite{McClenny2023} for a number of experiments. SA-PINNs weight each collocation point in the PINN loss function separately, and learn these weights as the PINN trains by combining gradient descent steps on the PINN's parameters with gradient ascent steps on the weights. SA-PINNs have been shown to improve the performance of PINNs when solving multi-scale problems, including the Helmholtz equation \cite{McClenny2023}. More specifically, implementing the boundary conditions as hard constraints, SA-PINNs use the following loss function to train the PINN,
\begin{equation} \label{eq:sapinn_loss_hard}
\mathcal{L}(\bm{\theta},\bm{\lambda}) =\frac{1}{N} \sum_{i=1}^{N} m(\lambda_i)(\mathcal{N}[\mathcal{C} u](\mathbf{x}_{i},\bm{\theta}) - f(\mathbf{x}_{i}))^{2}~,
\end{equation}
where $\{\lambda_i\}_{i=1}^{N}$ is a set of non-negative weights, with one weight for each collocation point, $m(\lambda_i)$ is a strictly increasing scalar function of $\lambda_i$, and the rest of the terms are as defined in~\cref{eq:pinn_loss_hard}. At each training step, gradient descent is used on~\cref{eq:sapinn_loss_hard} to update $\bm{\theta}$ and gradient ascent is used on~\cref{eq:sapinn_loss_hard} to update $\bm{\lambda}$. This results in the values $\lambda_i$ increasing during training and placing greater weight on collocation points which have consistently high PDE residuals.

In this work, for all tests using SA-PINNs we fix $m(\lambda_i)=\lambda_i$, initialize all $\lambda_i=1$ at the start of training, and use the same optimizer (Adam) and learning rate ($1\times10^{-3}$) to update $\bm{\theta}$ and $\bm{\lambda}$. For reference, the weights $\bm{\lambda}$ obtained after training the SA-PINNs used to solve the Helmholtz problems in \cref{fig:weak_helmholtz_constant_solution} are shown in \cref{fig:sapinn_weights}.

\begin{figure}[!t]
\centering
\includegraphics[width=0.94\textwidth]{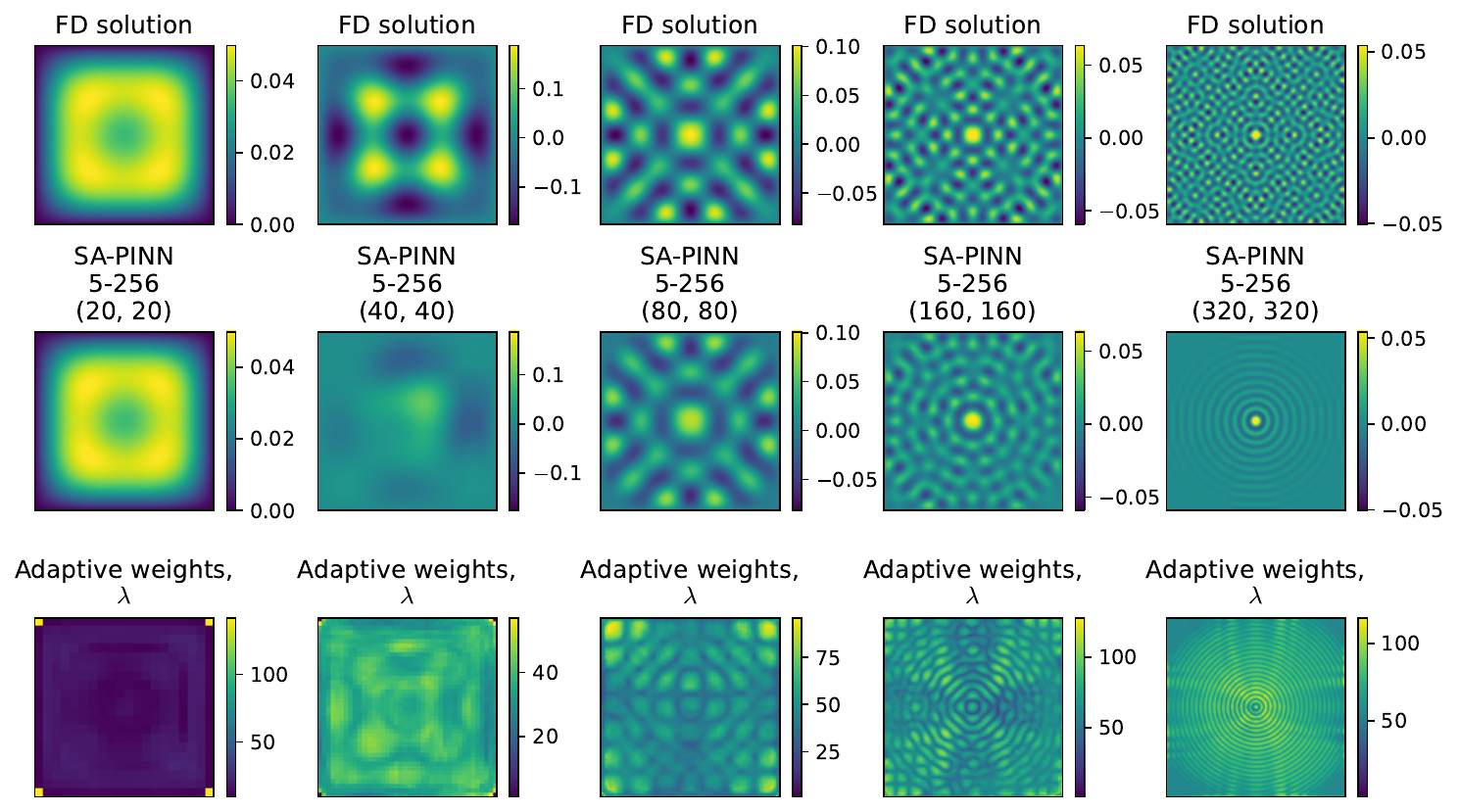}
\caption{\corr{FD reference solution (top row), SA-PINN solution (middle row), and SA-PINN weights $\bm{\lambda}$ (bottom row) obtained after training the SA-PINNs shown in \cref{fig:weak_helmholtz_constant_solution}.}
\label{fig:sapinn_weights}}
\end{figure}

\end{core}

\section{\corrr{FBPINN individual level contributions}}
\label{sec:levels}

\begin{figure}[!t]
\centering
\includegraphics[width=0.94\textwidth]{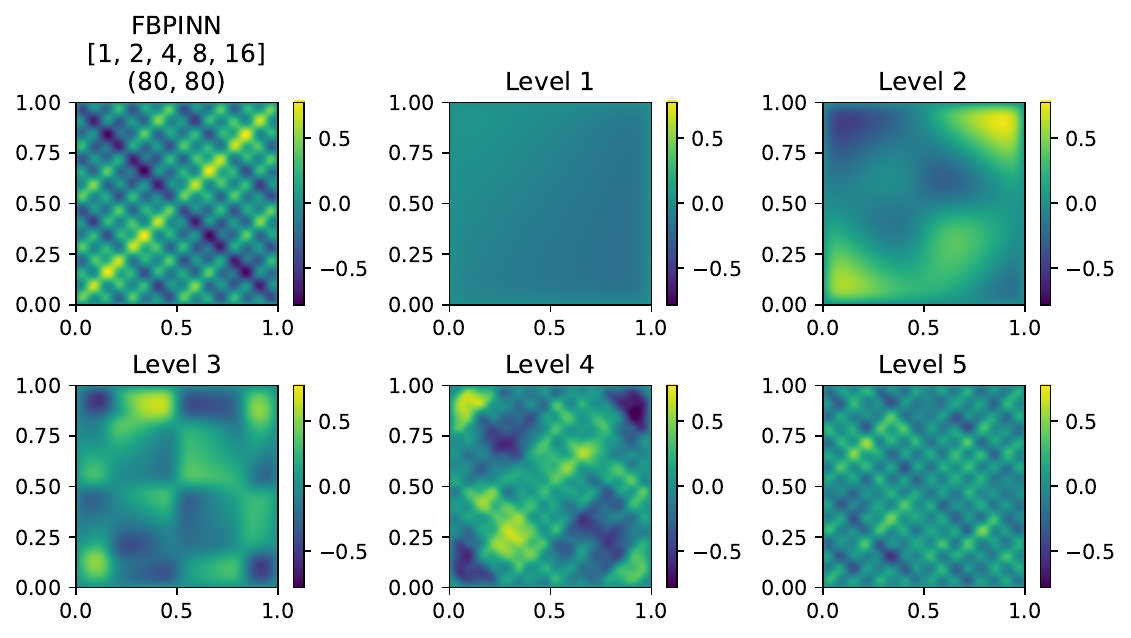}
\caption{\corrr{Contribution to the FBPINN solution from each level, for the 5-level FBPINN shown in \cref{fig:weak_laplace}. Top left shows the full FBPINN solution after summing the levels.}
\label{fig:levels}}
\end{figure}

\corrr{Figure~\ref{fig:levels} plots the contribution to the FBPINN solution from each level (the individual terms before summing across levels in \cref{eq:multilevel_fbpinn_network_architecture}), for the 5-level FBPINN shown in \cref{fig:weak_laplace}. We find that the finest level roughly learns the highest frequencies in the solution, whilst the coarser levels learn the lower frequency contributions. We expect this behaviour given the spectral bias of neural networks and the individual subdomain normalisation employed in the FBPINN.}

\bibliographystyle{elsarticle-num} 
\bibliography{references}







\end{document}